# CONTOUR REGRESSION: A GENERAL APPROACH TO DIMENSION REDUCTION[1]

By Bing Li, Hongyuan Zha and Francesca Chiaromonte

*Pennsylvania State University*

We propose a novel approach to sufficient dimension reduction in regression, based on estimating contour directions of small variation in the response. These directions span the orthogonal complement of the minimal space relevant for the regression and can be extracted according to two measures of variation in the response, leading to *simple* and *general contour regression* (SCR and GCR) methodology. In comparison with existing sufficient dimension reduction techniques, this contour-based methodology guarantees exhaustive estimation of the central subspace under ellipticity of the predictor distribution and mild additional assumptions, while maintaining $\sqrt{n}$-consistency and computational ease. Moreover, it proves robust to departures from ellipticity. We establish population properties for both SCR and GCR, and asymptotic properties for SCR. Simulations to compare performance with that of standard techniques such as ordinary least squares, sliced inverse regression, principal Hessian directions and sliced average variance estimation confirm the advantages anticipated by the theoretical analyses. We demonstrate the use of contour-based methods on a data set concerning soil evaporation.

**1. Introduction and background.** Consider the regression of a response $Y$ on a vector of continuous predictors $X = (X_1, \ldots, X_p)^T \in \mathbb{R}^p$. Sufficient dimension reduction is a body of theory and methods for reducing the dimension of $X$ while preserving information on the regression, that is, on the conditional distribution of $Y|X$ (see [7, 15, 16]). A *dimension reduction subspace* [3, 4] is defined as the column space of any $p \times d$ ($d \le p$) matrix $\eta$ such that

$$(1) \qquad\qquad Y \perp\!\!\!\perp X | \eta^T X,$$

Received July 2003; revised August 2004.

[1]Supported in part by NSF Grants DMS-02-04662 awarded to BL, DMS-03-11800 awarded to HZ and DMS-04-05681 awarded to FC, BL and HZ.

where $\perp\!\!\!\perp$ indicates independence. Thus, conditioning on $\eta^T X$, $Y$ and $X$ are independent, or equivalently, the conditional distribution of $Y|X$ equals that of $Y|\eta^T X$. Because relation (1) is unaffected by multiplying $\eta$ from the right by a nonsingular matrix, what matters is the column space of $\eta$ rather than its specific form. Also note that there can be many subspaces satisfying (1), because if it holds for $\eta$, then it also holds for any other matrix whose column space includes that of $\eta$. Naturally, we are interested in the subspace with the minimal dimension. Under mild conditions that are almost always verified in practice, the minimal subspace is uniquely defined and coincides with the intersection of all subspaces satisfying (1) (see [1, 4]). This intersection is called the *central subspace*, denoted by $\mathcal{S}_{Y|X}$, and its dimension is called the *structural dimension*, denoted by $q$.

The central subspace can be estimated without estimating a response surface, and without strong assumptions on the form of the dependence between $Y$ and $X$. Well-known estimation methods include ordinary least squares (OLS, [17]), sliced inverse regression (SIR, [15]; see also [9]), principal Hessian directions (PHD, [16]) and sliced average variance estimation (SAVE, [7]). These methods constitute effective premodeling tools to reduce high-dimensional regressions to equivalent ones comprising only a few linear combinations of the original predictors. Such a reduction greatly facilitates model building, as well as the use of nonparametric techniques. Dimension reduction methods also provide a comprehensive visualization of the data whenever the estimated structural dimension is 1, 2 or possibly 3, which is the case in a vast majority of practical applications. In this sense, sufficient dimension reduction provides a foundation for regression graphics, as argued in [4] and [1].

In many studies attention is restricted to the location component of the dependence between $Y$ and $X$, that is, to the regression function $E(Y|X)$. The *central mean subspace*, $\mathcal{S}_{E(Y|X)}$, was introduced by Cook and Li [5]. Because the conditional mean $E(Y|X)$ is determined by the distribution of $Y|X$, the central mean subspace is always contained in the central subspace. Cook and Li investigated the above mentioned methods in relation to their ability to estimate directions within the $\mathcal{S}_{E(Y|X)}$, and proposed alternative methods to target this subspace directly. See also [14].

The above methods enjoy the advantage of being computationally inexpensive and $\sqrt{n}$-consistent regardless of the original predictor dimension $p$ and the structural dimension $q$, thus avoiding the "curse of dimensionality" often affecting nonparametric techniques. The $\sqrt{n}$-consistency is achieved because these methods exploit *global* features of the dependence of $Y$ on $X$, in the sense they involve averaging among fixed portions of the data, regardless of the sample size. For instance, OLS employs sample moments, and SIR involves averaging predictors within slices of $Y$, where the size of each slice need not go to zero as $n \to \infty$.

The above methods also have common limitations, however. First, they require linear conditional means among predictors ([4], page 57), which is often imposed by assuming ellipticity of the distribution of $X$. When this condition fails, the estimators may converge to directions outside $\mathcal{S}_{Y|X}$. Second, even when linearity holds, OLS, PHD and SIR are not guaranteed to be exhaustive: they may converge at the $\sqrt{n}$-rate to a set of vectors that are in $\mathcal{S}_{Y|X}$ but do not span $\mathcal{S}_{Y|X}$. This lack of exhaustiveness is arguably one of the most important shortcomings of these methods. An instance is the heavy reliance of methods such as OLS and SIR on monotone trends in the dependence of $Y$ on $X$. For example, if $Y = (\beta^T X)^2 + \sigma\varepsilon$ with $\beta \in \mathbb{R}^p$ and $X \sim N(0, I_p)$, OLS and SIR will estimate 0 and therefore fail to detect $\beta$ itself. Based on early results obtained by Peters, Redner and Decell [19] in the special context of feature extraction, it can be shown that SAVE is indeed exhaustive if $X|Y$ is normally distributed. However, as we will see in Section 3, this assumption is very restrictive in the regression context. Thus, it is of both practical and theoretical significance to pursue exhaustive estimation under reasonably general sufficient conditions.

At the opposite end of the spectrum are adaptive methods that exploit *local* features of the dependence of $Y$ on $X$ [13, 23]. The strength of these methods is that they require much weaker assumptions (virtually none) on the distribution $X$. However, because they employ multivariate kernels that shrink with the sample size, their convergence rates are generally slower than $\sqrt{n}$. In addition, they are computationally intensive, as they iterate between nonparametric estimation of a multivariate unknown function and numerical maximization of the estimated function over a potentially high-dimensional space.

Here we propose a novel approach that targets contour directions, that is, directions along which the response surface is flat. Since contour directions span the orthogonal complement of the central subspace, estimating the former is equivalent to estimating the latter. We propose to extract contour directions according to two measures of variation in the response, leading to two methods: *simple* and *general contour regression* (SCR and GCR). Unlike traditional global methods such as OLS, SIR, PHD and SAVE, contour regression guarantees exhaustive estimation of the central subspace under ellipticity of $X$ and very mild additional assumptions. It also proves robust to violations of ellipticity. At the same time, unlike local methods, contour regression achieves $\sqrt{n}$-consistency regardless of the dimensions $p$ and $q$, and it is computationally inexpensive.

The remainder of the paper is organized as follows: Section 2 concerns population-level properties and asymptotic properties of SCR, and Section 3 presents sufficient conditions for exhaustiveness of SCR. Section 4 concerns population-level properties of GCR. Section 5 discusses the robustness of GCR against violations of predictor ellipticity. Section 6 presents simulations comparing the performance of SCR and GCR with that of OLS, SIR, PHD and SAVE. Section 7 reports the analysis of a data set. Section 8 contains final remarks. The main proofs are reported in the body of the paper, but some technical details are relegated to the Appendix.

## 2. Simple contour regression.

2.1. *Basic concepts.* Let $(X_1, Y_1), \ldots, (X_n, Y_n)$ be independent copies of the random pair $(X, Y)$, where $X \in \mathbb{R}^p$ and $Y \in \mathbb{R}$, let $F_{XY}$ be the joint distribution of $(X, Y)$, and let $F_n$ be the corresponding empirical distribution. We will be concerned with matrix-valued estimators of the form $T(F_n)$. If the columns of $T(F_{XY})$ belong to $\mathcal{S}_{Y|X}$, then we say that $T(F_n)$ is unbiased at the population level. If the columns of $T(F_{XY})$ actually span $\mathcal{S}_{Y|X}$, then we say that $T(F_n)$ is exhaustive at the population level. If $T(F_n)$ converges at the $\sqrt{n}$-rate to $T(F_{XY})$ in the first case, then we say that it is $\sqrt{n}$-consistent. If the $\sqrt{n}$-convergence holds in the second case, then we say that $T(F_n)$ is $\sqrt{n}$-exhaustive. In this section we will introduce simple contour regression and establish its exhaustiveness at the population level, as well as its $\sqrt{n}$-exhaustiveness.

To illustrate the basic intuition underlying contour regression, it is helpful to make an analogy between the *empirical directions* employed in this approach and *empirical distributions*. In defining an empirical distribution we put an equal probability mass at each observed point, based on the rationale that all the information about the random vector $X$ carried by the random sample $X_1, \ldots, X_n$ is captured by the positions of the observations themselves. With a similar rationale, it can be argued that all the *directional information* about $X$ carried by the data is captured by the empirical directions $(X_i - X_j)$, $i \neq j$, $i, j = 1, \ldots, n$. Since dimension reduction is about finding a specific set of directions—those along which the conditional distribution of $Y|X$ genuinely depends on $X$—it is natural to focus on these $\binom{n}{2}$ empirical directions. Roughly, contour regression extracts a subset of the empirical directions characterized by having small variation in $Y$, and then performs a principal component analysis on the extracted directions. Since contour directions form the orthogonal complement of the central subspace, an estimator of the latter is obtained from the components with smallest eigenvalues.

It is also important to remark that many global methods (see Introduction) gather directional information by slicing the response and processing predictor observations separately within each slice. This way, interslice directional information relevant to the regression, if any, is lost. On the other hand, empirical directions $(X_i - X_j)$ can "cut across" response slices, thus allowing contour regression to exploit interslice information.

For ease of exposition we will first introduce population-level quantities and then construct estimators by analogy.

2.2. *Population-level exhaustiveness.* Let $(\tilde{X}, \tilde{Y})$ be an independent copy of $(X, Y)$ and suppose that the central subspace $\mathcal{S}_{Y|X}$ for the regression of $Y$ on $X$ is spanned by the column space of a $p \times q$ matrix $\beta$ with $q < p$. Consider the matrix

$$K(c) = E[(\tilde{X} - X)(\tilde{X} - X)^T \mid |\tilde{Y} - Y| \le c].$$

We will show that the eigenvectors of $K(c)$ corresponding to its smallest $q$ eigenvalues span the central subspace. For this purpose we need the following assumption.

ASSUMPTION 2.1. For any choice of vectors $v \in \mathcal{S}_{Y|X}$ and $w \in (\mathcal{S}_{Y|X})^{\perp}$ such that $\|v\| = \|w\| = 1$, and some constant $c > 0$, we have

$$(2) \qquad \text{var}[w^T(\tilde{X} - X) \mid |\tilde{Y} - Y| \le c] > \text{var}[v^T(\tilde{X} - X) \mid |\tilde{Y} - Y| \le c].$$

This assumption is a reasonable one: because the conditional distribution of $Y|X$ depends on $v^T X$ but not on $w^T X$, we expect $Y$ to vary more with $v^T X$ than it does with $w^T X$. Hence, intuitively, within the same increment of $Y$, $w^T X$ should vary more than $v^T X$ does. In Section 3 we will prove that Assumption 2.1 holds under fairly general conditions. Here we give a few examples to illustrate its wide applicability. For reasons that will become clear later on, we will always impose this assumption on the standardized predictors with mean 0 and variance matrix $I_p$.

EXAMPLE 2.1. Suppose $X = (X_1, X_2)^T \sim N(0, I_2)$ and $Y = X_2^2 + \sigma \varepsilon$, with $\varepsilon \perp\!\!\!\perp X$ and $\varepsilon \sim N(0, 1)$. For this regression $\mathcal{S}_{Y|X}$ is the one-dimensional span of $\beta = (0, 1)^T$, and the conditional variances on the left- and right-hand sides of (2) are, respectively,

$$\lambda_1 = E[(\tilde{X}_1 - X_1)^2 \mid |\tilde{Y} - Y| \le c] \quad \text{and} \quad \lambda_2 = E[(\tilde{X}_2 - X_2)^2 \mid |\tilde{Y} - Y| \le c].$$

Because $X_1$ is independent of $X_2$ and $\varepsilon$, it is independent of $Y$. Therefore the first conditional expectation equals the unconditional expectation $E(X_1 - \tilde{X}_1)^2$, which is 2. We have computed the second conditional expectation numerically on a grid of values $c = 0.1, 0.5, 1, \ldots, 3$ and $\sigma = 0.1, 0.2, 0.3, \ldots, 2$. In all cases $\lambda_2 < \lambda_1 = 2$. Below is the tabulation corresponding to $\sigma = 0.3$ and selected values of $c$:

| $c$ | 0.10 | 0.50 | 1.00 | 1.50 | 2.00 | 2.50 | 3.00 |
|---|---|---|---|---|---|---|---|
| $\lambda_2$ | 0.85 | 0.90 | 1.04 | 1.20 | 1.34 | 1.45 | 1.54 |

We note that this is the case where OLS and SIR will estimate zero, providing no information about the central subspace.

EXAMPLE 2.2. Let $X = (X_1, X_2)^T$ and $\varepsilon$ be as in Example 2.1, and $Y = (X_2 - 1)^3 + \sigma\varepsilon$. We again used numerical integration to compute $\lambda_2$ for $c = 0.1, 0.5, 1, \ldots, 3$ and $\sigma = 0.1, 0.2, 0.3, \ldots, 2$, and again obtained values below 2 in all cases. A table for $\sigma = 0.3$ and selected values of $c$ follows:

| $c$ | 0.10 | 0.50 | 1.00 | 1.50 | 2.00 | 2.50 | 3.00 |
|-----|------|------|------|------|------|------|------|
| $\lambda_2$ | 0.24 | 0.28 | 0.36 | 0.45 | 0.53 | 0.61 | 0.67 |

We checked numerically numerous other response functions such as polynomials, exponential and logarithmic functions, trigonometric functions and so on, never encountering a violation of Assumption 2.1. In the next example, we verify Assumption 2.1 for a regression with a binary response. This example also helps to emphasize that Assumption 2.1 should be imposed on the standardized, rather than the original, predictor.

EXAMPLE 2.3. Let $Y \sim \text{Bernoulli}(1/2)$ and $X \in \mathbb{R}^2$. The regression is described through the inverse conditional distribution: $X|(Y = y) \sim N(\mu_y, I_2)$, with $\mu_0 = (0, -1)^T$ and $\mu_1 = (0, 1)^T$. Thus, here $\mathcal{S}_{Y|X}$ is again the one-dimensional span of $\beta = (0, 1)^T$ and the conditional variances on the left- and right-hand sides of Assumption 2.1 are $\lambda_1$ and $\lambda_2$ as defined in the previous examples. Simple calculations show that, for $0 \leq c \leq 1$,

$$E[(\tilde{X} - X)(\tilde{X} - X)^T | Y = \tilde{Y}] = E[(\tilde{X} - X)(\tilde{X} - X)^T | Y = \tilde{Y} = 0]$$
$$= 2E(XX^T | Y = 0) - 2E(X|Y = 0)E(X^T|Y = 0)$$
$$= 2\,\text{var}(X|Y = 0) = 2I_2.$$

Therefore, $\lambda_1 = \lambda_2 = 2$, and Assumption 2.1 fails. However, if we standardize the predictor vector at the outset, the assumption holds also for this binary regression. Note that $E(X) = 0$ and

$$\Sigma = \text{var}(X) = \begin{pmatrix} 1 & 0 \\ 0 & 2 \end{pmatrix}.$$

It follows that

$$E[(\tilde{Z} - Z)(\tilde{Z} - Z)^T | Y = \tilde{Y}]$$
$$= \Sigma^{-1/2} E[(\tilde{X} - X)(\tilde{X} - X)^T | Y = \tilde{Y}]\Sigma^{-1/2} = 2\Sigma^{-1}.$$

Because $\Sigma$ is diagonal, the central subspace $\mathcal{S}_{Y|Z}$ for the regression of $Y$ on the standardized predictor $Z = \Sigma^{-1/2}X$ is still the span of $(0, 1)^T$ (see, e.g., [4], page 106). Thus, on the $Z$-scale $\lambda_1 = 1 < 2 = \lambda_2$.

In the next theorem and corollary we prove that if $X$ is elliptically contoured and Assumption 2.1 holds, then the population vectors from SCR exhaust the central subspace $\mathcal{S}_{Y|X}$. We first consider the standardized $X$.

THEOREM 2.1. *Suppose that $X$ has an elliptical distribution with $E(X) = 0$ and $\mathrm{var}(X) = I_p$. If Assumption 2.1 holds, then the eigenvectors of $K(c)$ corresponding to its smallest $q$ eigenvalues span the central subspace $\mathcal{S}_{Y|X}$.*

PROOF. Suppose, without loss of generality, that the $q$ columns of $\beta$ form an orthonormal set in $\mathbb{R}^p$. Let $\gamma_1, \ldots, \gamma_{p-q}$ be an orthonormal basis for $(\mathcal{S}_{Y|X})^\perp$. We need to show that (i) $\gamma_1, \ldots, \gamma_{p-q}$ are eigenvectors of $K(c)$, and (ii) their corresponding eigenvalues are the largest among its eigenvalues. If (i) and (ii) hold, then the eigenvectors of $K(c)$ corresponding to its smallest $q$ eigenvalues will coincide with $\mathrm{span}(\beta) = \mathcal{S}_{Y|X}$, as desired.

(i) Let $B = (\gamma_1, \ldots, \gamma_{p-q}, \beta_1, \ldots, \beta_q)$ and $W = B^T X$. By Proposition 6.3 of [4], page 106, the central subspace for the regression of $Y$ on $W$, $\mathcal{S}_{Y|W}$, is spanned by the vectors $B^{-1}\beta_1, \ldots, B^{-1}\beta_q$. Note that $B^{-1}\beta_i = e_{p-q+i}$, that is, the vector in $\mathbb{R}^p$ with a 1 in the $(p-q+i)$th position and 0's in all remaining ones. By construction $W$ still has a spherical distribution with mean 0 and variance $I_p$. Next, let $(\tilde{X}, \tilde{Y})$ be an independent copy of $(X, Y)$, $\tilde{W} = B^T \tilde{X}$, and

$$K_1(c) = E[(\tilde{W} - W)(\tilde{W} - W)^T \mid |\tilde{Y} - Y| \le c].$$

It is easy to see that $K_1(c) = B^{-1}K(c)B$, and hence that $\gamma_i$ is an eigenvector of $K(c)$ if and only if $B^T\gamma_i$ is an eigenvector of $K_1(c)$. Thus, it suffices to show that $e_1, \ldots, e_{p-q}$ are eigenvectors of $K_1(c)$.

Let $\phi_i : \mathbb{R}^p \mapsto \mathbb{R}^p$, $i = 1, \ldots, p-q$, be the map that changes the sign of the $i$th element of a vector in $\mathbb{R}^p$ [e.g., $\phi_1(x) = (-x_1, x_2, \ldots, x_p)^T$]. Observe that $\phi_i$ is an invertible mapping with $\phi_i^{-1} = \phi_i$, and that $x, \phi_1(x), \ldots, \phi_{p-q}(x)$ are on the same sphere centered at the origin. Let $C$ be any measurable set. Then

$$\Pr(W \in C) = \Pr(W \in \phi_i(C)) = \Pr(\phi_i^{-1}(W) \in C) = \Pr(\phi_i(W) \in C),$$

where the first equality follows because $W$ has a spherical distribution and $x \mapsto \phi_i(x)$ is an orthogonal transformation, and the third equality follows because $\phi_i = \phi_i^{-1}$. Thus $W$ and $\phi_i(W)$ have the same distribution. Furthermore, because of the conditional independence,

$$Y \perp\!\!\!\perp (W_1, \ldots, W_p)|W_{p-q+1}, \ldots, W_p,$$

the conditional distributions of $Y|W$ and $Y|\phi_i(W)$ (where $i = 1, \ldots, p-q$) are both equal to the conditional distribution of $Y|W_{p-q+1}, \ldots, W_p$. It follows that $(W, Y)$ and $(\phi_i(W), Y)$ have the same distribution. Consequently,

$$\begin{aligned}
K_1(c) &= E[(\tilde{W} - W)(\tilde{W} - W)^T \mid |\tilde{Y} - Y| \le c] \\
&= E[(\phi_i(\tilde{W}) - \phi_i(W))(\phi_i(\tilde{W}) - \phi_i(W))^T \mid |\tilde{Y} - Y| \le c].
\end{aligned}$$

So the matrix $K_1(c)$ can be re-expressed as the average of the right-hand sides of the first and the second lines above; that is,

$$E[\tfrac{1}{2}(\tilde{W} - W)(\tilde{W} - W)^T + \tfrac{1}{2}(\phi_i(\tilde{W}) - \phi_i(W))(\phi_i(\tilde{W}) - \phi_i(W))^T \mid |\tilde{Y} - Y| \le c].$$

Denote the above matrix by $E(A_i \mid |\tilde{Y} - Y| \le c)$ in the obvious way. Then it is easy to see that the $(i, i)$th element of $A_i$ is $(\tilde{W}_i - W_i)^2$ and the $(i, j)$th and $(j, i)$th elements are identically 0 whenever $j \ne i$. Hence, for $i = 1, \ldots, p - q$, $e_i$ is an eigenvector of $K_1(c)$ corresponding to the eigenvalue $\lambda_i = E((W_i - \tilde{W}_i)^2 \mid |Y - \tilde{Y}| \le c)$.

(ii) From (i) it follows that $\gamma_1, \ldots, \gamma_{p-q}$ are eigenvectors of $K(c)$ corresponding to the eigenvalues $\lambda_1, \ldots, \lambda_{p-q}$. Let $\gamma_{p-q+1}, \ldots, \gamma_p$ be the other $q$ eigenvectors of $K(c)$, corresponding to the eigenvalues $\lambda_{p-q+1}, \ldots, \lambda_p$. By construction these $q$ vectors span $\mathcal{S}_{Y|X}$. The eigenvalues $\lambda_i$, $i = 1, \ldots, p$, can be expressed as

$$\lambda_i = \gamma_i^T K(c) \gamma_i = \mathrm{var}[\gamma_i^T(X - \tilde{X}) \mid |Y - \tilde{Y}| \le c].$$

But because $\gamma_1, \ldots, \gamma_{p-q}$ are vectors in $(\mathcal{S}_{Y|X})^{\perp}$ and $\gamma_{p-q+1}, \ldots, \gamma_p$ are vectors in $\mathcal{S}_{Y|X}$, Assumption 2.1 suffices to ensure that $\lambda_1, \ldots, \lambda_{p-q}$ are the largest eigenvalues of $K(c)$. □

Let us now turn to $X$ with an arbitrary elliptically contoured distribution having mean $\mu$ and nonsingular variance matrix $\Sigma$. Theorem 2.1 asserts that if we let $Z = \Sigma^{-1/2}(X - \mu)$, then the eigenvectors $\gamma_{p-q+1}, \ldots, \gamma_p$ corresponding to the smallest $q$ eigenvalues of the matrix $E[(\tilde{Z} - Z)(\tilde{Z} - Z)^T \mid |\tilde{Y} - Y| \le c]$ span $\mathcal{S}_{Y|Z}$. Consequently, by Proposition 6.3 of [4], page 106, the vectors $\Sigma^{-1/2}\gamma_{p-q+1}, \ldots, \Sigma^{-1/2}\gamma_p$ span $\mathcal{S}_{Y|X}$. Thus we have the following generalization of Theorem 2.1:

COROLLARY 2.1. *Suppose that $X$ has an elliptically contoured distribution with mean $\mu$ and nonsingular variance matrix $\Sigma$. Suppose that Assumption 2.1 holds with $X$ and $\tilde{X}$ replaced by $\Sigma^{-1/2}(X - \mu)$ and $\Sigma^{-1/2}(\tilde{X} - \mu)$. Let $\gamma_{p-q+1}, \ldots, \gamma_p$ be the eigenvectors of the matrix*

$$\Sigma^{-1/2} E[(\tilde{X} - X)(\tilde{X} - X)^T \mid |\tilde{Y} - Y| \le c] \Sigma^{-1/2}$$

*corresponding to its smallest $q$ eigenvalues. Then the vectors $\Sigma^{-1/2}\gamma_{p-q+1}, \ldots, \Sigma^{-1/2}\gamma_p$ span $\mathcal{S}_{Y|X}$.*

Note that since the corollary postulates Assumption 2.1 on the standardized predictor, it does also apply to regressions with discrete responses, such as the one described in Example 2.3.

2.3. *Estimation and $\sqrt{n}$-exhaustiveness.* We now construct a sample estimate of the matrix $K(c)$. As before let $(\tilde{X}, \tilde{Y})$ be an independent copy of $(X, Y)$, and consider the matrix

$$(3) \qquad H(c) = E[(\tilde{X} - X)(\tilde{X} - X)^T I(|\tilde{Y} - Y| \leq c)].$$

Since $K(c)$ and $H(c)$ differ only by the proportionality constant $\Pr(|\tilde{Y} - Y| \leq c)$, their eigenvectors coincide. We will thus consider an estimate of $H(c)$ for simplicity. Let $(X_1, Y_1), \ldots, (X_n, Y_n)$ be an independent sample from the random pair $(X, Y)$. The estimating procedure will mimic the theoretical development in Section 2.2:

(a) Compute sample mean and variance matrix of the predictor $X$:

$$\hat{\mu} = n^{-1} \sum_{i=1}^{n} X_i, \qquad \widehat{\Sigma} = n^{-1} \sum_{i=1}^{n} (X_i - \hat{\mu})(X_i - \hat{\mu})^T.$$

(b) Compute the matrix-valued $U$-statistic:

$$(4) \qquad \widehat{H}(c) = \frac{1}{\binom{n}{c}} \sum_{(i,j) \in N} (X_j - X_i)(X_j - X_i)^T I(|Y_j - Y_i| \leq c),$$

where $N$ is the index set $\{(i, j) : i = 2, \ldots, n; j = 1, \ldots, i - 1\}$.

(c) Compute the spectral decomposition of $\widehat{\Sigma}^{-1/2} \widehat{H}(c) \widehat{\Sigma}^{-1/2}$ and let $\hat{\gamma}_{p+q-1}, \ldots, \hat{\gamma}_p$ be the eigenvectors corresponding to the smallest $q$ eigenvalues.

(d) The span of these eigenvectors estimates $\mathcal{S}_{Y|Z}$, where $Z$ is the standardized version of $X$. Thus, our estimate of the central subspace is

$$\widehat{\mathcal{S}}_{Y|X} = \mathrm{span}(\widehat{\Sigma}^{-1/2} \hat{\gamma}_{p-q+1}, \ldots, \widehat{\Sigma}^{-1/2} \hat{\gamma}_p).$$

In practice, as for other sufficient dimension reduction methods, a testing procedure will be needed to determine statistically how many eigenvectors to take. Here we assume the dimension $q$ of the central subspace to be known, leaving the development of an appropriate testing procedure for future investigation. However, in the next section we lay the groundwork for constructing such a procedure by further exploring the eigenvalue structure of the matrix $K(c)$. Before doing so, we demonstrate the $\sqrt{n}$-exhaustiveness of SCR.

THEOREM 2.2. *Suppose that $\Sigma$ is nonsingular and that the components of $X$ have finite fourth moments. Then*

$$\widehat{\Sigma}^{-1/2} \widehat{H}(c) \widehat{\Sigma}^{-1/2} = \Sigma^{-1/2} H(c) \Sigma^{-1/2} + O_p(n^{-1/2}).$$

PROOF. Since $X$ has finite fourth moment, $\widehat{\Sigma}$ is a $\sqrt{n}$-consistent estimator of $\Sigma$ by the central limit theorem. Since $\Sigma$ is nonsingular, we have that $\widehat{\Sigma}^{-1/2}$ is a $\sqrt{n}$-consistent estimator of $\Sigma^{-1/2}$ by the continuous mapping theorem. It follows that

$$\widehat{\Sigma}^{-1/2}\widehat{H}(c)\widehat{\Sigma}^{-1/2} - \Sigma^{-1/2}H(c)\Sigma^{-1/2}$$
$$= \Sigma^{-1/2}(\widehat{H}(c) - H(c))\Sigma^{-1/2} + O_p(n^{-1/2}).$$

Next, let $\mathrm{vec}(\cdot)$ be the operator that stacks the columns of a matrix [for $A$ with columns $a_1, \ldots, a_k$, $\mathrm{vec}(A) = (a_1^T, \ldots, a_k^T)^T$], and let $\mathrm{vec}^T(\cdot)$ denote the transpose of $\mathrm{vec}(\cdot)$. Note that $\widehat{H}(c)$ is a matrix-valued $U$-statistic with elements having finite second moments. By the central limit theorem for $U$ statistics ([20], Chapter 5) it easily follows that $\sqrt{n}\,\mathrm{vec}(\widehat{H}(c) - H(c))$ converges in distribution to a $p^2$-dimensional multivariate normal vector with mean 0 and variance matrix

$$4E[\mathrm{vec}(H_1(X, Y; c))\mathrm{vec}^T(H_1(X, Y; c))] - 4\,\mathrm{vec}(H(c))\mathrm{vec}^T(H(c)),$$

where $H_1(X, Y; c)$ is the conditional expectation

$$E[(\tilde{X} - X)(\tilde{X} - X)^T I(|\tilde{Y} - Y| \le c)|X, Y].$$

Consequently $\widehat{H}(c) - H(c) = O_p(1/\sqrt{n})$, which completes the proof.  □

As a consequence of this theorem, $\hat{\gamma}_{p-q+1}, \ldots, \hat{\gamma}_p$ provide a $\sqrt{n}$-exhaustive estimator of $\mathcal{S}_{Y|Z}$, and hence $\widehat{\Sigma}^{-1/2}\hat{\gamma}_{p-q+1}, \ldots, \widehat{\Sigma}^{-1/2}\hat{\gamma}_p$ provide a $\sqrt{n}$-exhaustive estimator of $\mathcal{S}_{Y|X}$.

The role played by the constant $c$ here is similar to the width of a slice in SIR and SAVE. It differs from the width of a kernel in a typical nonparametric estimator in that it need not go to 0 as $n \to \infty$. The thresholding can actually be implemented in two ways: fixing a numerical value for $c$, or fixing a proportion of empirical directions $X_i - X_j$ [out of $\binom{n}{2}$]. This distinction is relevant for theoretical analysis and simulation. We find that using the 5–15% empirical directions ranking lowest in terms of response absolute difference works well in simulation studies (see Section 6). A more careful investigation of the thresholding rule is important, but goes beyond the scope of the present article—we expect good thresholding to depend on the dimension $p$, $q$, and possibly other factors.

The asymptotic analyses we present here are all carried out for the thresholding based on a fixed value of $c$. However, they can be easily paralleled for thresholding based on a fixed proportion: Modulo the fact that $c$ need not go to 0 as $n \to \infty$, the comparison of these two thresholding options is analogous to that between a kernel and a nearest-neighbor estimator [18, 21, 22].

2.4. *Toward testing hypothesis.* In this section we will show that, under an additional assumption, the largest $p - q$ eigenvalues of $\Sigma^{-1/2}K(c)\Sigma^{-1/2}$ are identically 2, so that $\gamma_1, \ldots, \gamma_{p-q}$ are the eigenvectors of $2I_p - \Sigma^{-1/2}K(c)\Sigma^{-1/2}$ corresponding to eigenvalues equal to 0. This paves the way for constructing a test statistic to determine the dimension of $\mathcal{S}_{Y|Z}$ (and hence $\mathcal{S}_{Y|X}$), because the problem now is converted into testing how many of the smallest eigenvalues of $2I_p - \Sigma^{-1/2}K(c)\Sigma^{-1/2}$ are 0. Tests of this type can be constructed using the asymptotic distribution of small singular values developed by Eaton and Tyler [10]. Similar tests have been constructed in other contexts in [15, 16]; [4], Chapter 11 and [2, 5, 14]. We expect that a test statistic and related sampling distribution for SCR can be obtained analogously. The additional assumption is usually referred to as the constant conditional variance assumption, and is often evoked when developing such tests:

ASSUMPTION 2.2. If $\beta$ is a matrix whose columns form a basis in $\mathcal{S}_{Y|X}$, then the conditional variance $\text{var}(X|\beta^T X)$ is a nonrandom matrix.

The next lemma states various implications of conditional independence, which will be used in the subsequent development.

LEMMA 2.1. (a) *If* $V_1, \ldots, V_6$ *are random vectors satisfying* $(V_1, V_2, V_3) \perp\!\!\!\perp (V_4, V_5, V_6)$, $V_1 \perp\!\!\!\perp V_2|V_3$ *and* $V_4 \perp\!\!\!\perp V_5|V_6$, *then* $(V_1, V_4) \perp\!\!\!\perp (V_2, V_5)|(V_3, V_6)$.

(b) *If* $V_1, \ldots, V_4$ *are random vectors satisfying* $(V_1, V_2) \perp\!\!\!\perp (V_3, V_4)$, *then* $V_1 \perp\!\!\!\perp V_3|(V_2, V_4)$.

(c) *If* $V_1, V_2, V_3$ *are random vectors satisfying* $(V_1, V_2) \perp\!\!\!\perp V_3$, *then* $V_1 \perp\!\!\!\perp V_3|V_2$.

Part (c) is a special case of a well-known result; see, for example, [8] and [4], Proposition 4.6. The proofs of (a) and (b) are similar to those used in these papers, and are given in the Appendix.

THEOREM 2.3. *Suppose that* $X$ *has an elliptical distribution and that Assumptions* 2.1 *and* 2.2 *hold. Let* $q$ *be the dimension of* $\mathcal{S}_{Y|X}$. *Then the* $p - q$ *largest eigenvalues of* $\Sigma^{-1/2}K(c)\Sigma^{-1/2}$ *are identically* 2.

PROOF. Let $\beta$ be a $p \times q$ matrix whose columns span $\mathcal{S}_{Y|X}$, and $Z = \Sigma^{-1/2}(X - \mu)$. Then the columns of $\eta = \Sigma^{1/2}\beta$ span $\mathcal{S}_{Y|Z}$ and $\Sigma^{-1/2}K(c)\Sigma^{-1/2}$ is the matrix

$$K_1(c) = E[(\tilde{Z} - Z)(\tilde{Z} - Z)^T \mid |\tilde{Y} - Y| \le c] \equiv E(\Delta_1 \Delta_1^T \mid |\Delta_2| \le c),$$

where we have abbreviated $\tilde{Z} - Z$ and $\tilde{Y} - Y$ by $\Delta_1$ and $\Delta_2$, respectively. Because $(Z, Y, \eta^T Z) \perp\!\!\!\perp (\tilde{Z}, \tilde{Y}, \eta^T \tilde{Z})$, $Z \perp\!\!\!\perp Y | \eta^T Z$ and $\tilde{Z} \perp\!\!\!\perp \tilde{Y} | \eta^T \tilde{Z}$, we have, by Lemma 2.1(a), that $(Z, \tilde{Z}) \perp\!\!\!\perp (Y, \tilde{Y}) | \eta^T Z, \eta^T \tilde{Z}$. This in turn implies that $\Delta_1 \perp\!\!\!\perp \Delta_2 | \eta^T Z, \eta^T \tilde{Z}$. Consequently,

$$E(\Delta_1 \Delta_1^T \mid |\Delta_2| \le c) = E[E(\Delta_1 \Delta_1^T | \eta^T Z, \eta^T \tilde{Z}, |\Delta_2| \le c) \mid |\Delta_2| \le c]$$
$$= E[E(\Delta_1 \Delta_1^T | \eta^T Z, \eta^T \tilde{Z}) \mid |\Delta_2| \le c].$$

The conditional expectation inside the brackets on the right-hand side can be decomposed as the sum of four terms:

$$(5) \qquad \begin{aligned} & E(ZZ^T | \eta^T Z, \eta^T \tilde{Z}) - E(\tilde{Z} Z^T | \eta^T Z, \eta^T \tilde{Z}) \\ & \quad - E(Z \tilde{Z}^T | \eta^T Z, \eta^T \tilde{Z}) + E(\tilde{Z} \tilde{Z}^T | \eta^T Z, \eta^T \tilde{Z}). \end{aligned}$$

Since $(Z, \eta^T Z) \perp\!\!\!\perp \eta^T \tilde{Z}$, we have $Z \perp\!\!\!\perp \eta^T \tilde{Z} | \eta^T Z$ by Lemma 2.1(c). Hence the first term becomes

$$E(ZZ^T | \eta^T Z) = \mathrm{var}(Z | \eta^T Z) + E(Z | \eta^T Z) E(Z^T | \eta^T Z).$$

Let $P = \eta(\eta^T \eta)^{-1} \eta^T$ be the orthogonal projection onto $\mathcal{S}_{Y|Z}$ and let $Q = I - P$ be the orthogonal projection onto $(\mathcal{S}_{Y|Z})^\perp$. Then it can be shown by Assumption 2.2 that $\mathrm{var}(Z | \eta^T Z) = Q$. Because $Z$ has a spherical distribution, $E(Z | \eta^T Z) = PZ$, so that the first term in (5) reduces to $Q + PZZ^T P$. By Lemma 2.1(b), the second term in (5) factorizes into $-E(\tilde{Z} | \eta^T Z, \eta^T \tilde{Z}) E(Z^T | \eta^T Z, \eta^T \tilde{Z})$, which by Lemma 2.1(c) further reduces to $-E(\tilde{Z} | \eta^T \tilde{Z}) E(Z^T | \eta^T Z)$. Hence, again using sphericity of $Z$, the second term in (5) is $-P\tilde{Z} Z^T P$. By similar arguments the third and fourth terms in (5) are $-PZ\tilde{Z}^T P$ and $Q + P\tilde{Z}\tilde{Z}^T P$, respectively. Therefore

$$(6) \qquad E(\Delta_1 \Delta_1^T \mid |\Delta_2| \le c) = 2Q + PE(\Delta_1 \Delta_1^T \mid |\Delta_2| \le c)P.$$

Let $v$ be a vector in $(\mathcal{S}_{Y|Z})^\perp$, and multiply this matrix by $v^T$ from the left and by $v$ from the right to obtain $\mathrm{var}(v^T \Delta_1 \mid |\Delta_2| \le c) = 2$. This completes the proof.  □

Note that without Assumption 2.1 the theorem still holds to an extent, in the sense that the eigenvectors of $E(\Delta_1 \Delta_1^T \mid |\Delta_2| \le c)$ orthogonal to span$(\eta)$ still have eigenvalues equal to 2. However, without this assumption exhaustiveness would be lost, because we cannot rule out the possibility that eigenvalues other than these may also be 2.

This eigenvalue structure is similar to that of SAVE. For some $c > 0$ and $y$, let $S(c)$ be the sliced averaged variance $S(c) = \mathrm{var}(X \mid |Y - y| \le c)$. In this notation $I_p - S(c)$ is the SAVE matrix for a slice centered at $y$. Under ellipticity and Assumption 2.2, we have

$$S(c) = Q + PS(c)P$$

(see [7]). Thus the eigenvalues corresponding to the eigenvectors of $S(c)$ that are orthogonal to the central subspace are identically 1. However, there are important differences between contour regression methods and SAVE; these will be briefly discussed at the end of Section 4.3.

**3. Sufficient conditions for exhaustive estimation.** In order to place the theory of simple contour regression on a firmer foundation, we devote this section to deriving a sufficient condition for Assumption 2.1. As shown in the previous sections, if this assumption holds, then SCR provides $\sqrt{n}$-*exhaustive* estimation of the central subspace $\mathcal{S}_{Y|X}$; that is, the estimating vectors converge with $\sqrt{n}$-rate to a set of vectors that span $\mathcal{S}_{Y|X}$ in its entirety. Sufficient conditions of this type are extremely elusive; to our knowledge none has been established with reasonable generality for other $\sqrt{n}$-consistent methods such as OLS, PHD, SIR or SAVE. Results from an early, prescient paper by Peters, Redner and Decell [19] lead to exhaustiveness of SAVE under the condition that $X|Y$ is multivariate normal. However, this condition is very restrictive—note that even in a typical location regression of the form $Y = f(X) + \varepsilon$ with $X$ and $\varepsilon$ independent and both normally distributed, this assumption is not met unless $f(\cdot)$ is linear. Because of its generality, the sufficient condition given here for SCR is the first of its kind.

We will need the notion of stochastic ordering. Let $S$ and $T$ be two random variables. We say that $S$ is stochastically less than or equal to $T$ if, for any real number $r$, $\Pr(S \leq r) \geq \Pr(T \leq r)$, and write this as $S \leq_d T$. If, in addition, the inequality is strict on a subset of the real line with positive Lebesgue measure, we say that $S$ is stochastically (strictly) less than $T$ and write $S <_d T$. The following lemma is obvious, and its proof will be omitted.

LEMMA 3.1.   *Suppose that $S$ and $T$ are random variables taking values in a common set $\Omega \subset \mathbb{R}$, and that $S <_d T$. Then:*

(a) $E(S) < E(T)$.
(b) *Given a monotone real-valued function $g : \Omega \mapsto \mathbb{R}$, $g(S) <_d g(T)$ if $g(\cdot)$ is increasing, and $g(T) <_d g(S)$ if $g(\cdot)$ is decreasing.*

As a special case, consider a pair of random variables $(S, T)$. We will write $(S|T = t_1) <_d (S|T = t_2)$ if, for any $r$, $\Pr(S \leq r|T = t_1) \geq \Pr(S \leq r|T = t_2)$, with strict inequality held on a set with positive Lebesgue measure.

The following lemma, which is proved in the Appendix, will also be used.

LEMMA 3.2.   *Let $p(s)$ and $q(s)$ be the densities of nonnegative random variables $S$ and $T$ taking values in a common support $\Omega \subset \mathbb{R}^+$, and suppose that $p(s)/q(s)$ is decreasing in $s$. Then $E(S) < E(T)$.*

In developing a sufficient condition for Assumption 2.1, we restrict ourselves to a *location structure*, that is, to regressions of the kind

$$(7) \qquad Y = f(\beta^T X) + \sigma \varepsilon, \qquad \varepsilon \perp\!\!\!\perp X, E(\varepsilon) = 0.$$

Ultimately, the sufficient condition will be imposed on the behavior of $f(\cdot)$. Let $(\tilde{X}, \tilde{\varepsilon})$ be an independent copy of $(X, \varepsilon)$, $\Delta = \tilde{X} - X$, $T = \tilde{\varepsilon} - \varepsilon$, and let $F_T(\cdot)$ be the cumulative distribution function of $T$. For the statement of the following theorem, it will be more informative to write $f(\beta^T x)$ merely as $g(x)$.

THEOREM 3.1. *Suppose that $X$ has an elliptically contoured distribution with $E(X) = 0$ and $\mathrm{var}(X) = I_p$, and that Assumption 2.2 holds. Moreover, suppose that model (7) holds with the density $f_T(t)$ of $F_T(t)$ being a decreasing function of $|t|$. If for any $\alpha \in \mathcal{S}_{Y|X}$ and whenever $0 \le \delta_1 < \delta_2$ we have*

$$(8) \quad |g(X + \Delta) - g(X)| \,|\, \{|\alpha^T \Delta| = \delta_1\} <_d |g(X + \Delta) - g(X)| \,|\, \{|\alpha^T \Delta| = \delta_2\},$$

*then Assumption 2.1 holds for every $c > 0$.*

Before proving the theorem, let us comment on its significance. To understand the intuition behind condition (8), first consider the case where $X$ is a scalar random variable. Intuitively, condition (8) should hold trivially if $g$ is a monotone function, because it holds pointwise in $X = x$ with $<_d$ replaced by ordinary inequality $<$ (see Example 3.1 below). However, condition (8) by no means restricts $g(\cdot)$ to being monotone, because being stochastically large or small is an average behavior for all values of $X$, and does not require being large or small for every single value $X = x$. It then does seem to make sense to assume that $g(X + \Delta)$ is collectively farther away from $g(X)$ if $\Delta$ is larger: this is simply requiring $g$ to be reasonably variable. In the multivariate case, condition (8) requires this to hold along any direction $\alpha$ in the space $\mathcal{S}_{Y|X}$, which is the space along which $g(x)$ does vary. Also the requirement that $f_T(t)$ decreases with $|t|$ is not a severe restriction, considering that this density is symmetric about 0 by construction. Finally, the result can be generalized straightforwardly to nonstandardized elliptical predictors. Thus, Theorem 3.1 allows us to conclude that, for elliptical predictors with constant conditional variance along the central subspace (Assumption 2.2), all that is required to guarantee Assumption 2.1—and hence exhaustiveness of SCR—are very mild conditions on the behavior of the mean function and the error term in (7) (some instances are provided after the proof).

PROOF OF THEOREM 3.1. Let $\alpha \in \mathcal{S}_{Y|X}$ and $\xi \in (\mathcal{S}_{Y|X})^\perp$. Theorem 2.3 implies that $\mathrm{var}(\xi^T(\tilde{X} - X) \,|\, |\tilde{Y} - Y| \le c) = 2$, which is the same as the unconditional variance $\mathrm{var}(\alpha^T(\tilde{X} - X))$. Hence, it suffices to show that

$$\mathrm{var}(\alpha^T(\tilde{X} - X) \,|\, |\tilde{Y} - Y| \le c) < \mathrm{var}(\alpha^T(\tilde{X} - X)).$$

Let $U = |\tilde{Y} - Y|$ and $V = (\alpha^T(\tilde{X} - X))^2$. We are then to show that $E(V|U \leq c) < E(V)$. Let $f_V(\cdot)$ be the density of $V$. Then

$$(9) \qquad E(V|U \leq c) = \int_0^\infty v \frac{\Pr(U \leq c|V = v)}{\Pr(U \leq c)} f_V(v)\, dv.$$

Now, let $r(v) = \Pr(U \leq c|V = v)/\Pr(U \leq c)$. Then $r(v)f_V(v)$ is itself a density on $\mathbb{R}^+$. By Lemma 3.2, if we can show that $r(v)$ is a decreasing function of $v$, then the right-hand side of (9) is smaller than $\int v f_V(v)\, dv$ and the proof is complete. So let us show that $\Pr(U \leq c|V = v)$ decreases in $v$. Note that

$$\Pr(U \leq c|V = v) = E\{\Pr(U \leq c|X, \tilde{X})|V = v\}.$$

Because $(X, \tilde{X})$ is independent of $(\varepsilon, \tilde{\varepsilon})$ we can re-express the conditional probability $\Pr(U \leq c|X, \tilde{X})$ as

$$\Pr(g(X) - g(\tilde{X}) - c \leq T \leq g(X) - g(\tilde{X}) + c)$$
$$= F_T(g(X) - g(\tilde{X}) + c) - F_T(g(X) - g(\tilde{X}) - c).$$

Because the roles of $X$ and $\tilde{X}$ can be exchanged, and because $V = (\alpha^T(\tilde{X} - X))^2 = (\alpha^T(X - \tilde{X}))^2$, we have

$$\Pr(U \leq c|V = v) = E[F_T(g(X) - g(\tilde{X}) + c) - F_T(g(X) - g(\tilde{X}) - c)|V = v]$$
$$= E[F_T(g(\tilde{X}) - g(X) + c) - F_T(g(\tilde{X}) - g(X) - c)|V = v].$$

So $\Pr(U \leq c|V = v)$ can be written as the average of the two expressions on the right-hand sides of the first and second equalities in the above display. That is, if we write $g(X) - g(\tilde{X}) + c$ as $A$ and $g(X) - g(\tilde{X}) - c$ as $-B$, then $\Pr(U \leq c|V = v)$ can be written as $E(F_T(A) - F_T(-B) + F_T(B) - F_T(-A))/2$. However, since $F_T(\cdot)$ is the cumulative distribution function of a symmetric density, we have $F_T(t) - F_T(-t) = 2F_T(t) - 1$ for any $t$. Hence

$$\Pr(U \leq c|V = v)$$
$$= E[F_T(g(X) - g(\tilde{X}) + c) + F_T(g(\tilde{X}) - g(X) + c)|V = v] - 1$$
$$= E[F_T(R + c) + F_T(-R + c)|V = v] - 1 \equiv E(G(R)|V = v) - 1,$$

where $R = |g(\tilde{X}) - g(X)|$. Thus it suffices to show that $E(G(R)|V = v)$ decreases with $v$. Because $f_T(\cdot)$ is symmetric about 0,

$$G'(r) = f_T(r + c) - f_T(-r + c) = f_T(|r + c|) - f_T(|r - c|).$$

But since $r$ and $c$ are both positive, $r + c > |r - c|$. It follows that the right-hand side is negative, and $G(r)$ is strictly decreasing in $r$ for $r \geq 0$. By Lemma 3.1(a), $E(G(R)|V = v)$ will be decreasing in $v$ if, for any $v_2 > v_1 \geq 0$,

$$(10) \qquad (G(R)|V = v_1) <_d (G(R)|V = v_2).$$

However, because $G(R)$ is a decreasing function of $R$, by Lemma 3.1(b) inequality (10) will hold if $(R|V = v_1) <_d (R|V = v_2)$. The latter inequality is equivalent to (8). $\quad\square$

To illustrate the generality of this sufficient condition we now verify it for some examples. Let us consider the following specialization of our location structure. In (7) take $X \sim N(0, I_2)$ and $\beta = (0, 1)^T$, so that $Y = f(X_2) + \sigma\varepsilon$. Consider the conditional probability

$$\Pr(|f(\tilde{X}_2) - f(X_2)| \le r \,|\, |\tilde{X}_2 - X_2| = \delta).$$

Condition (8) will be satisfied if, for each $r > 0$, this quantity decreases in $\delta$. Because the distribution of $\tilde{X}_2 - X_2$ is symmetric about zero, this probability is

$$\Pr(|f(\tilde{X}_2) - f(X_2)| \le r | \tilde{X}_2 - X_2 = \delta)/2$$
$$+ \Pr(|f(\tilde{X}_2) - f(X_2)| \le r | \tilde{X}_2 - X_2 = -\delta)/2.$$

Because the roles of $X$ and $\tilde{X}$ can be interchanged, the conditioning argument $\tilde{X}_2 - X_2 = -\delta$ in the second term can be replaced by $\tilde{X}_2 - X_2 = \delta$, and hence

$$\Pr(|f(\tilde{X}_2) - f(X_2)| \le r \,|\, |\tilde{X}_2 - X_2| = \delta)$$
$$= \Pr(|f(\tilde{X}_2) - f(X_2)| \le r | \tilde{X}_2 - X_2 = \delta) \equiv \Psi(\delta).$$

Thus (8) will hold if, for each $r > 0$, $\Psi(\delta)$ is a decreasing function of $\delta > 0$. Now $X_2$ and $\tilde{X}_2$ can be written as $S + T$ and $S - T$ where $S = (X_2 + \tilde{X}_2)/2$ and $T = (X_2 - \tilde{X}_2)/2$. Note that, by normality of $X$, $S$ and $T$ are independent. Hence,

$$\Psi(\delta) = \Pr(|f(S+T) - f(S-T)| \le r \,|\, T = \delta) = \Pr(|f(S+\delta) - f(S-\delta)| \le r).$$

So for this specialization of our location structure we only need to verify

$$|f(S + \delta_1) - f(S - \delta_1)| <_d |f(S + \delta_2) - f(S - \delta_2)|$$

for $S \sim N(0, 1/2)$ and for any $0 \le \delta_1 < \delta_2$. The following examples both reduce to verifying this inequality.

EXAMPLE 3.1.   Suppose that $f(x_2)$ is a continuous and monotone function which without loss of generality can be assumed to be monotone increasing. Then for any $\delta_1 < \delta_2$, $|f(s + \delta_1) - f(s - \delta_1)| \le_d |f(s + \delta_2) - f(s - \delta_2)|$. To see that the strict inequality $(<_d)$ holds, let $r$ be a real number in the set $\{f(s + \delta_1) - f(s - \delta_1) : s \in \mathbb{R}\}$. By continuity there is an $s_0$ such that $f(s_0 + \delta_1) - f(s_0 - \delta_1) = r < f(s_0 + \delta_2) - f(s_0 - \delta_2)$. Hence, in the neighborhood of $(s_0 - \tau, s_0)$, $f(s + \delta_1) - f(s - \delta_1) < r < f(s + \delta_2) - f(s - \delta_2)$. Writing $|f(S + \delta) - f(S - \delta)|$ as $R(\delta)$, we have

$$\Pr(R(\delta_1) \leq r) = \Pr(R(\delta_1) \leq r, R(\delta_2) \leq r) + \Pr(R(\delta_1) \leq r, R(\delta_2) > r)$$
$$= \Pr(R(\delta_2) \leq r) + \Pr(R(\delta_1) \leq r, R(\delta_2) > r).$$

Because the set $\{s : R(\delta_1) \leq r, R(\delta_2) > r\}$ contains an open interval, $\Pr(R(\delta_1) \leq r, R(\delta_2) > r) > 0$, and consequently $\Pr(R(\delta_1) \leq r) < \Pr(R(\delta_2) \leq r)$. Because $\Pr(R(\delta_1) \leq r) - \Pr(R(\delta_2) \leq r)$ is continuous in $r$, this inequality holds in an open interval around $r$, which has positive Lebegue measure.

EXAMPLE 3.2. Let $f(x_2) = (x_2 - a)^2$. Example 2.1 is a special case of this regression with $a = 0$ and $\varepsilon \sim N(0, \sigma^2)$. In this case $R(s; \delta) = 4|s - a|\delta$. Hence for $0 \leq \delta_1 < \delta_2$, $R(s; \delta_1) < R(s; \delta_2)$ for all $s$. By an argument similar to the one in Example 3.1, it is easy to see that $R(S; \delta_1) <_d R(S; \delta_2)$.

## 4. General contour regression.

### 4.1. *Estimation.*  The idea underlying SCR is to use the inequality $|Y - \tilde{Y}| \leq c$ to identify vectors aligned with the contour directions. However, this inequality also picks up other directions when the regression function is nonmonotone. Under ellipticity such directions are averaged out, so that the method remains $\sqrt{n}$-exhaustive. Nevertheless, these "wrong" directions do tend to decrease efficiency by blurring up the "right" ones. In other words, the inequality $|Y - \tilde{Y}| \leq c$ is not a very sensitive contour identifier for nonmonotone functions—even though it is sufficiently sensitive to maintain $\sqrt{n}$-exhaustiveness. We now illustrate this point using the regression in Example 2.1.

To construct the left panel of Figure 1, we generated twenty observations $(X_i, Y_i)$, $i = 1, \ldots, 20$, according to the regression in Example 2.1, with $\sigma = 0.3$. We then used the threshold value $c = 0.5$, connecting by a solid line segment any two points $X_i, X_j \in \mathbb{R}^2$ satisfying $|Y_i - Y_j| \leq 0.5$. Roughly speaking (note that we have ignored the rescaling issue which has little bearing on this discussion), SCR picks up the contour directions by a principal component analysis of the vectors represented by these line segments. We see that, though most of the segments are horizontal (i.e., aligned with the true contour direction), there are a considerable number of segments pointing to arbitrary directions. This is because the response surface is $U$-shaped and the inequality $|Y_i - Y_j| \leq 0.5$ does not discriminate between the segments aligned with the contour and those across the $U$-shaped surface that also have small increments in $Y$. Though the arbitrary directions tend to average out due to the ellipticity of the distribution of $X$, they make the picture less sharp and the method less efficient.

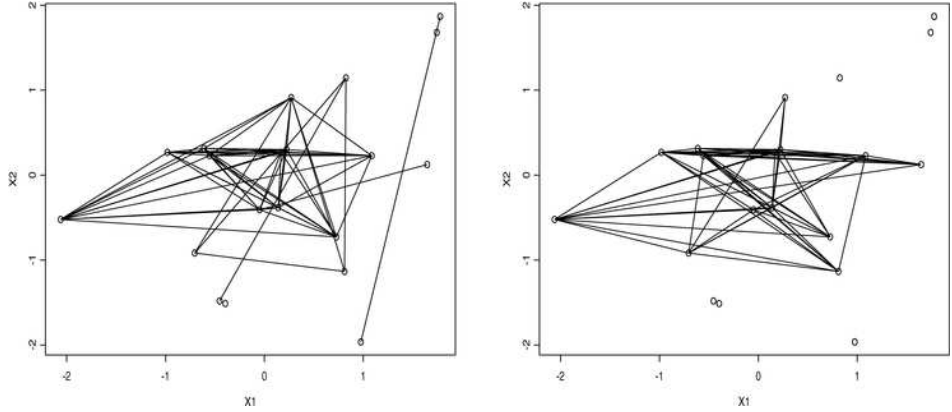

FIG. 1. *Directions identified by* $|Y - \tilde{Y}| \leq c$ (left panel) *and those identified by* $\widehat{V}(X_i, X_j, \rho) \leq c$ (right panel).

To overcome this drawback we replace the contour identifier $|Y_i - Y_j| \leq c$ by a more sensitive one. Consider the variance of $Y$ along the line through $x_i$ and $x_j$. Formally, let $\ell(t; x_i, x_j) = (1 - t)x_i + tx_j$, $t \in \mathbb{R}$, be the straight line that goes through $x_i$ and $x_j$, and define

$$V(x_i, x_j) = \text{var}(Y | X = \ell(t; x_i, x_j) \text{ for some } t).$$

For a more concrete expression, let $\delta(x_i, x_j)$ be the $p \times (p - 1)$ matrix $(\delta_1, \ldots, \delta_{p-1})$ whose columns form a basis in $(x_j - x_i)^\perp$. Then $V(x_i, x_j)$ can be re-expressed as

$$(11) \qquad V(x_i, x_j) = \text{var}(Y | \delta^T(x_i, x_j)X = \delta^T(x_i, x_j)x_i).$$

We will aim at identifying contour vectors by the smallness of this conditional variance.

The next task is to construct a sample estimate of $V(X_i, X_j)$. We will denote the line $\ell(\cdot; X_i, X_j)$ by $\ell(X_i, X_j)$. For any $X_k$, let $d(X_k, \ell(X_i, X_j))$ be the Euclidean distance between $X_k$ and the line $\ell(X_i, X_j)$; that is,

$$d(X_k, \ell(X_i, X_j)) = \min_{t \in \mathbb{R}} \|X_k - \ell(t; X_i, X_j)\|,$$

where $\|\cdot\|$ stands for the Euclidean norm. Because $\|X_k - \ell(t; X_i, X_j)\|^2$ is a quadratic function of $t$, this minimum distance can be expressed explicitly as

$$d(X_k, \ell(X_i, X_j)) = \left[ \|X_k - X_i\|^2 - \frac{\{(X_k - X_i)^T(X_j - X_i)\}^2}{\|X_j - X_i\|^2} \right]^{1/2}.$$

For any $\rho > 0$, we define the tube of radius $\rho$ connecting $X_i$ and $X_j$ to be the set

$$C_{ij}(\rho) = \{X_k : d(X_k, \ell(X_i, X_j)) \leq \rho, k = 1, \ldots, n\}.$$

According to this definition, each tube contains at least two points in the sample. Next we estimate the variance of $Y$ along these tubes. Let $n_{ij}(\rho)$ be the number of points in the tube $C_{ij}(\rho)$, and let

$$\widehat{V}(X_i, X_j; \rho) = \frac{1}{n_{ij}(\rho)} \sum_{X_k \in C_{ij}(\rho)} (Y_k - \bar{Y}_{ij}(\rho))^2,$$

$$\text{where } \bar{Y}_{ij}(\rho) = \frac{1}{n_{ij}(\rho)} \sum_{X_k \in C_{ij}(\rho)} Y_k.$$

We can now identify the contour directions by the smallness of $\widehat{V}(X_i, X_j; \rho)$.

Plotted in the right panel of Figure 1 are the same sample points as in the left panel, but with the line segments picked up by $\widehat{V}(X_i, X_j; \rho) \leq c$, where $c = 0.5$ and $\rho = 0.3$. We can see that many of the segments pointing to random directions in the left panel have been removed. To get a quantitative comparison, we calculated the first principal component for the line segments in each panel, which equals $(0.9169, 0.3991)^T$ for the left panel and $(0.9991, -0.0417)^T$ for the right panel. The latter is much closer to the direction $(1, 0)^T$, the population vector orthogonal to $\mathcal{S}_{Y|X}$.

We now construct the estimator of $\mathcal{S}_{Y|X}$. Along lines similar to those followed in Section 2, we standardize the predictor observations to $\widehat{Z}_i = \widehat{\Sigma}^{-1/2}(X_i - \hat{\mu})$, and form the matrix

$$(12) \qquad \widehat{F}(c) = \frac{1}{\binom{n}{2}} \sum_{(i,j) \in N} (\widehat{Z}_j - \widehat{Z}_i)(\widehat{Z}_j - \widehat{Z}_i)^T I(\widehat{V}(\widehat{Z}_i, \widehat{Z}_j; \rho) \leq c),$$

where $N$ is the same index set as used in (4). The matrix $\widehat{F}(c)$ takes the place of $\widehat{\Sigma}^{-1/2}\widehat{H}(c)\widehat{\Sigma}^{-1/2}$ for the simple contour regression. As in SCR, we take the spectral decomposition of $\widehat{F}(c)$, and use $\hat{\gamma}_{p+q-1}, \ldots, \hat{\gamma}_p$, the eigenvectors corresponding to the smallest $q$ eigenvalues, to form

$$\hat{\mathcal{S}}_{Y|X} = \text{span}(\widehat{\Sigma}^{-1/2}\hat{\gamma}_{p-q+1}, \ldots, \widehat{\Sigma}^{-1/2}\hat{\gamma}_p).$$

Regarding the choice of $c$, comments similar to those made at the end of Section 2.3 apply here. In particular, as a rule of thumb we propose to use 5% to 15% of the $\binom{n}{2}$ empirical directions.

4.2. *Population-level exhaustiveness.* Assume that $X$ is already standardized to $E(X) = 0$ and $\text{var}(X) = I_p$ (so $Z$ is $X$ itself). The population version of the matrix $\widehat{F}(c)$ in (12) is

$$F(c) = E[(X - \tilde{X})(X - \tilde{X})^T I(V(X, \tilde{X}) \leq c)],$$

which is proportional to the matrix

$$G(c) = E[(X - \tilde{X})(X - \tilde{X})^T | V(X, \tilde{X}) \leq c].$$

Here we will demonstrate that, for sufficiently small $c$, the eigenvectors corresponding to the smallest $q$ eigenvalues of $G(c)$ span $\mathcal{S}_{Y|X}$. For this purpose we introduce an assumption that parallels Assumption 2.1. Again $(\tilde{X}, \tilde{Y})$ indicates an independent copy of $(X, Y)$.

ASSUMPTION 4.1.  For any choice of vectors $v \in \mathcal{S}_{Y|X}$ and $w \in (\mathcal{S}_{Y|X})^{\perp}$ such that $\|v\| = \|w\| = 1$, and some constant $c > 0$, we have

$$(13) \quad \operatorname{var}[w^T(\tilde{X} - X)|V(X, \tilde{X}) \le c] > \operatorname{var}[v^T(\tilde{X} - X)|V(X, \tilde{X}) \le c].$$

The interpretation of this assumption is similar to that of Assumption 2.1, except that $V(X, \tilde{X})$ replaces $|\tilde{Y} - Y|$ as the measure of variation of $Y$ along the line through $X$ and $\tilde{X}$. We now deduce population exhaustiveness under this assumption. Once again we do so for a spherical predictor without loss of generality.

THEOREM 4.1.  *Suppose that $X$ has an elliptical distribution with $E(X) = 0$ and $\operatorname{var}(X) = I_p$. Then, under Assumption 4.1, the eigenvectors of $G(c)$ corresponding to its smallest $q$ eigenvalues span the central subspace $\mathcal{S}_{Y|X}$.*

The proof of this theorem is similar to that of Theorem 2.1 and will be given in the Appendix. The generalization to an arbitrary elliptical distribution for $X$ is similar to Corollary 2.1 and will be omitted.

We expect $\sqrt{n}$-consistency to hold also for GCR estimation. However, the asymptotic analysis is substantially more complex than for the SCR case, because the estimator cannot be rendered directly as a $U$-statistic. Alternative techniques must be developed for such an analysis. In this paper we will not prove $\sqrt{n}$-convergence rate for GCR, but will back up our claim by simulation in Section 6.

4.3. *Sufficient conditions for exhaustive estimation.*  Next, following a reasoning similar to that in Section 3, and again in reference to the location structure in (7), we derive a sufficient condition for Assumption 4.1.

Note that, since $\mathcal{S}_{Y|X} = \operatorname{span}(\beta)$, for a $p \times r$ matrix $\delta$ $(r \le p)$ we will have

$$(14) \quad \operatorname{var}(f(\beta^T X)|\delta^T X) > 0$$

unless $\operatorname{span}(\beta) \subset \operatorname{span}(\delta)$; that is, $f(\beta^T X)$ is not a function of $\delta^T X$ unless $\delta$ spans a space containing the central subspace.

THEOREM 4.2.  *Suppose that $X$ has an elliptically contoured distribution with $E(X) = 0$ and $\operatorname{var}(X) = I_p$, and that model (7) holds. Then Assumption 4.1 is satisfied for all sufficiently small $c > 0$ for which $\{(x, \tilde{x}) : V(x, \tilde{x}) \le c\}$ is a nonempty set.*

Proof. We first show that $V(x, \tilde{x}) \geq \sigma^2$ for all $x$ and $\tilde{x}$ and that equality holds for some $x$ and $\tilde{x}$, so that, whenever $c \geq \sigma^2$, $\{(x, \tilde{x}) : V(x, \tilde{x}) \leq c\}$ is nonempty. Let $\delta$ be any $p \times r$ matrix with $r \leq p$. Because $\varepsilon \perp\!\!\!\perp (\beta^T X, \delta^T X)$, we have by Lemma 2.1(c) $\beta^T X \perp\!\!\!\perp \varepsilon | \delta^T X$. Hence

$$(15) \qquad \begin{aligned} \operatorname{var}(Y | \delta^T X = t) &= \operatorname{var}(f(\beta^T X) | \delta^T X = t) + \operatorname{var}(\varepsilon | \delta^T X = t) \\ &= \operatorname{var}(f(\beta^T X) | \delta^T X = t) + \sigma^2, \end{aligned}$$

where for the second equality we have used the independence between $\varepsilon$ and $X$. Now take $\delta = \delta(x, \tilde{x})$ and $t = \delta^T(x, \tilde{x})x$, where $\delta(x, \tilde{x})$ is as defined above display (11). We see that $V(x, \tilde{x}) \geq \sigma^2$, and that equality holds whenever $(\tilde{x}-x)$ is orthogonal to $\operatorname{span}(\beta) = \mathcal{S}_{Y|X}$. With this in mind the assertion of the theorem can be rewritten as: if $v \in \mathcal{S}_{Y|X}$ and $w \in (\mathcal{S}_{Y|X})^\perp$, then for sufficiently small $\tau > 0$

$$\operatorname{var}[w^T(\tilde{X} - X) | V(X, \tilde{X}) \leq \sigma^2 + \tau] > \operatorname{var}[v^T(\tilde{X} - X) | V(X, \tilde{X}) \leq \sigma^2 + \tau].$$

By the definition of conditional expectation,

$$\lim_{\tau \downarrow 0} \operatorname{var}[(\tilde{X} - X) | V(X, \tilde{X}) \leq \sigma^2 + \tau] = \operatorname{var}[(\tilde{X} - X) | V(X, \tilde{X}) = \sigma^2].$$

Hence if we can show that

$$(16) \quad \operatorname{var}[v^T(\tilde{X} - X) | V(X, \tilde{X}) = \sigma^2] < \operatorname{var}[w^T(\tilde{X} - X) | V(X, \tilde{X}) = \sigma^2],$$

then the inequality will hold for all sufficiently small $\tau > 0$, proving the theorem.

To prove (16), note that (15) also implies that $\operatorname{var}(Y | \delta^T X = t) = \sigma^2$ if and only if $\operatorname{var}(f(\beta^T X) | \delta^T X = t) = 0$. However, because of (14), this will happen if and only if $\operatorname{span}(\beta) \subset \operatorname{span}(\delta)$. Taking $\delta = \delta(x, \tilde{x})$ and $t = \delta^T(x, \tilde{x})x$, we see that $V(x, \tilde{x}) = \sigma^2$ if and only if $\operatorname{span}(\beta) \subset \operatorname{span}(\delta(x, \tilde{x}))$, which is equivalent to $(x - \tilde{x}) \perp \operatorname{span}(\beta)$. Hence the left-hand side of (16) equals 0.

It remains to show that the right-hand side of (16) is positive. First, note that the roles of $X$ and $\tilde{X}$ are exchangeable, and hence

$$E[(X - \tilde{X}) | V(X, \tilde{X}) = \sigma^2] = E[(\tilde{X} - X) | V(\tilde{X}, X) = \sigma^2].$$

However, by the definition of $V(x, \tilde{x})$, $V(x, \tilde{x}) = V(\tilde{x}, x)$ for all $x$ and $\tilde{x}$, and hence

$$E[(X - \tilde{X}) | V(X, \tilde{X}) = \sigma^2] = E[(\tilde{X} - X) | V(X, \tilde{X}) = \sigma^2].$$

It follows that both sides of this equation must be 0, and so the right-hand side of (16) reduces to $E[(w^T(\tilde{X} - X))^2 | V(X, \tilde{X}) = \sigma^2]$. If this quantity were 0, then $w^T(x - \tilde{x}) = 0$ whenever $V(x, \tilde{x}) = \sigma^2$, which holds whenever $\beta^T(x - \tilde{x}) = 0$. Because the support of $X$ is spherical, $(x - \tilde{x})$ can run through every direction in $\mathbb{R}^p$. Hence $\operatorname{span}(\beta)^\perp \subset \operatorname{span}(w)^\perp$, or equivalently $w \in \operatorname{span}(\beta)$, which is a contradiction. $\quad\square$

The conditions in Theorem 4.2 are much weaker than those in Theorem 3.1. Predictor ellipticity and the structure in (7) are postulated in both cases. However, constant conditional variance (Assumption 2.2) is not required in Theorem 4.2, and essentially no requirement is posed on the behavior of the mean function and the error term in (7). Thus, GCR will be exhaustive under settings even more general than those required by SCR. Intuitively, this is because $V(X, \tilde{X})$ possesses stronger discriminating power than $|\tilde{Y} - Y|$: it can identify the contour vectors of any function $f(\beta^T X)$, as long as the latter genuinely depends on all the components of $\beta^T X$ [which is the case because of the minimality of the central subspace $S_{Y|X} = \text{span}(\beta)$ discussed in Section 1].

As mentioned at the end of Section 2.4, we now briefly discuss differences between contour regression methods and SAVE. First, in Section 3 and this section we have shown that SCR and GCR are guaranteed to be exhaustive under population-level conditions much milder than the ones assumed for exhaustiveness of SAVE. Second, contour regression methods break the barriers of slices, making more efficient use of data. In comparison, slice-based methods such as SAVE cannot exploit interslice information. Third, whereas the construction of SCR is somewhat similar to that of SAVE, and we suspect the gain in accuracy of SCR (which will be demonstrated by simulation) to be largely due to its efficient use of interslice information, GCR differs more intrinsically from slice-based methods: it employs a more sensitive contour identifier, and is thereby capable of picking up directions not easily detected by SAVE or SCR when the regression surface is complex. In Section 6 we show by simulation how this leads to improved accuracy in estimating the central subspace.

**5. Robustness against nonellipticity.** The population exhaustiveness of our contour-based methodology relies on ellipticity of the predictor distribution. This is because in the theoretical development we have treated the constant $c$ in (4) and (12) as fixed with respect to the sample size $n$. Ellipticity of the distribution of $X$ helps to balance out the effect of line segments that are not aligned with the contour directions. As mentioned in the Introduction, ellipticity requirements are ubiquitous for global methods such as OLS, SIR, PHD and SAVE. They are adopted to guarantee linear relationships among predictors, which in turn are needed for the methods to estimate directions within the central subspace. When the number of predictors $p$ is relatively small, diagnosing and remedying departures from ellipticity is relatively straightforward—in practice, scatterplot matrices are used to search for marked curvatures, and predictor transformations or data reweighting to mitigate such curvatures [4, 6]. However, especially when $p$ is large, diagnosing and remedying departures from ellipticity becomes laborious and complicated.

Notwithstanding the theoretical requirement, contour regression methods (especially GCR, whose contour identifier is more sensitive) can perform well even under violations of ellipticity. In Section 6 we will address this robustness by simulation; here we motivate it from a theoretical viewpoint. We will show that, postulating again the location structure in (7), the eigenvectors corresponding to the smallest $p - q$ eigenvalues of the matrix

$$A = E[(\tilde{X} - X)(\tilde{X} - X)^T | V(X, \tilde{X}) = \sigma^2]$$

span the orthogonal complement of the central subspace, $(\mathcal{S}_{Y|X})^\perp$, even when $X$ is not elliptical. This suggests that if we let $c$ decrease to $\sigma^2$ as $n$ increases, then the eigenvectors corresponding to the smallest $p - q$ eigenvalues of $\hat{F}(c)$ in (12) (after appropriate transformation by $\hat{\Sigma}^{-1/2}$) will tend to recover the whole $\mathcal{S}_{Y|X}$, regardless of the shape of the distribution of $X$. In practice, if we make $c$ small [i.e., close to the smallest value of $\hat{V}(\hat{Z}_i, \hat{Z}_j; \rho)$ in (12)], then GCR is likely to estimate the central subspace exhaustively and effectively even if the shape of $X$ does not help the process by averaging out erroneous directions, as is the case under ellipticity.

THEOREM 5.1. *Suppose that model (7) holds and that $X$ is a continuous random vector with an open support $\mathcal{X} \subset \mathbb{R}^p$. Then the matrix $A$ has exactly $p - q$ zero eigenvalues, and their corresponding eigenvectors span $(\mathcal{S}_{Y|X})^\perp$. In symbols,*

$$\ker(A) = \mathcal{S}_{Y|X},$$

*where $\ker(A) = \{h \in \mathbb{R}^p : Ah = 0\}$ is the kernel of $A$.*

PROOF. Note that $(\tilde{X} - X)$ is orthogonal to $\text{span}(\beta) = \mathcal{S}_{Y|X}$ if and only if $\text{span}(\beta) \subset \text{span}(\delta(X, \tilde{X}))$, which, by the argument following (15), happens if and only if $V(X, \tilde{X}) = \sigma^2$. Thus, conditioning on $V(X, \tilde{X}) = \sigma^2$, $(\tilde{X} - X)$ is orthogonal to $\text{span}(\beta)$. It follows that, whenever $h$ belongs to $\text{span}(\beta)$, $Ah = 0$, and thus $\text{span}(\beta) \subset \ker(A)$.

Conversely, suppose $h$ belongs to $\ker(A)$. Then

$$h^T A h = E[(h^T(\tilde{X} - X))^2 | V(X, \tilde{X}) = \sigma^2] = 0.$$

Thus, whenever $V(X, \tilde{X}) = \sigma^2$, $h$ is orthogonal to $(\tilde{X} - X)$. Equivalently, whenever $(\tilde{X} - X)$ is orthogonal to $\text{span}(\beta)$, $(\tilde{X} - X)$ is orthogonal to $h$. In other words, if we let $\mathcal{X}^* = \{\tilde{x} - x : \tilde{x} \in \mathcal{X}, x \in \mathcal{X}\}$, then

$$\mathcal{X}^* \cap (\text{span}(\beta))^\perp \subset \mathcal{X}^* \cap (\text{span}(h))^\perp.$$

However, because $\mathcal{X}$ is an open set, $\mathcal{X}^*$ is an open set containing 0. By Lemma A.1 in the Appendix $(\text{span}(\beta))^\perp \subset (\text{span}(h))^\perp$, or equivalently $h \in \text{span}(\beta)$, as desired. □

Intuitively, the theorem shows that the *only* directions $(x - \tilde{x})$ along which the variance $V(x, \tilde{x})$ achieves its minimum are those aligned with the contour. This is largely due to conditioning on the conditional variance, a population quantity. An analogous result cannot be derived for SCR.

Theorem 5.1 also suggests that, when we are not confident about the ellipticity of the distribution of $X$, we should use a stricter thresholding in the analysis [i.e., choose a small value of $c$, or include a small proportion of the $\binom{n}{2}$ empirical directions]. This makes contour regression estimators more similar to kernel estimators, whose consistency depends on the kernel width approaching 0 as $n \to \infty$. We will return to this point in the Conclusions.

**6. Simulation results.** We now compare the performance of both versions of contour regression, SCR and GCR, with that of well-known dimension reduction methods ensuring $\sqrt{n}$-consistency, such as OLS, SIR, PHD and SAVE. For such comparisons, we need to introduce a measure of distance between two subspaces of $\mathbb{R}^p$. Let $\mathcal{S}_1$ and $\mathcal{S}_2$ be two $q$-dimensional subspaces of $\mathbb{R}^p$ and let $P_{\mathcal{S}_1}$, $P_{\mathcal{S}_2}$ be the orthogonal projections onto $\mathcal{S}_1$ and $\mathcal{S}_2$, respectively. We use the distance

$$\text{dist}(\mathcal{S}_1, \mathcal{S}_2) = \|P_{\mathcal{S}_1} - P_{\mathcal{S}_2}\|,$$

where $\|\cdot\|$ is the Euclidean norm, that is, the maximum singular value of a matrix.

In the following, we present five examples covering a range of possible regression contexts. For both SCR and GCR we need to determine the number of empirical directions to include in the principal component analysis, and for GCR we also need to determine the tube radius $\rho$. Though in this paper we will not deal with the optimal choice of these numbers, related issues will be discussed to some extent in the examples.

In the first three examples the sample size $n$ and dimension $p$ are relatively small, whereas in the last two examples they are much larger. Instead of using a fixed $c$ for thresholding, we fix the proportion $r$ of the number of empirical directions with smallest variation (absolute response differences for SCR or tube variances for GCR) relative to $\binom{n}{2}$, the total number of empirical directions. For the first three examples we use $r = 6qn / \binom{n}{2}$ for SCR and $r = 2qn / \binom{n}{2}$ for GCR, and use $\rho = 1$ for GCR. For the last two examples we use $r = 5\%$ for both SCR and GCR and $\rho = 2$ for GCR.

EXAMPLE 6.1. Consider the regression

$$(17) \qquad\qquad Y = X_1^2 + X_2 + \sigma\varepsilon,$$

where $X \sim N(0, I_4)$, so that predictor ellipticity holds, $\varepsilon \sim N(0,1)$ and $\varepsilon \perp\!\!\!\perp X$. Here the central subspace is of dimension $q = 2$ and is spanned by the vectors $(1,0,0,0)^T$ and $(0,1,0,0)^T$. We compare SCR and GCR with SIR, SAVE and PHD using three different values of the error standard deviations: $\sigma = 0.1$, $0.4$ and $0.8$. Because OLS can pick up at most one direction, it is not included in this comparison. For each value of $\sigma$ we draw 500 samples of size $n = 100$, and on each sample we apply the five methods to produce five estimates of $\mathcal{S}_{Y|X}$. Next we compute the distance between these estimates and the true central subspace according to the definition at the beginning of this section. Finally, we compute an average and a standard error from the resulting 500 distances, for each $\sigma$ value and estimation method. Results are presented in Table 1 (DIST and SE columns correspond to average and standard error of the distances, resp.).

The numbers in Table 1 indicate that both SCR and GCR outperform SIR, SAVE and PHD in this example. Intuitively this is because SIR does not perform well when there is no linear trend, thus failing to pick up the second direction $(0,1,0,0)^T$, whereas PHD, and to a lesser extent SAVE, do not perform well when there is no quadratic trend, thus failing to give accurate estimates of the first direction $(1,0,0,0)^T$. In contrast, both SCR and GCR, as also demonstrated theoretically, provide comprehensive estimates of the central subspace. Note that SAVE performs better than SIR and PHD—by inspecting a few typical cases (results not presented) we find that SAVE does a better job at picking up the linear trend than PHD. Nevertheless, it is much less accurate than SCR and GCR. From the table we can also see that GCR generally outperforms SCR.

In the next example predictor ellipticity is maintained, but the comparison is based on a more complex regression surface in which linear and quadratic trends are not neatly separated along two coordinate directions. In this more complicated case, SIR, SAVE and PHD can also detect both directions.

EXAMPLE 6.2. Consider the regression

$$Y = X_1/(0.5 + (X_2 + 1.5)^2) + (1 + X_2)^2 + \sigma\varepsilon,$$

TABLE 1
*Comparison of SCR, GCR and other methods for Example 6.1*

| | SCR | | GCR | | SIR | | SAVE | | PHD | |
|---|---|---|---|---|---|---|---|---|---|---|
| $\sigma$ | DIST | SE | DIST | SE | DIST | SE | DIST | SE | DIST | SE |
| 0.1 | 0.23 | 0.11 | 0.16 | 0.07 | 0.78 | 0.24 | 0.43 | 0.25 | 0.80 | 0.21 |
| 0.4 | 0.25 | 0.11 | 0.20 | 0.08 | 0.79 | 0.23 | 0.54 | 0.27 | 0.79 | 0.21 |
| 0.8 | 0.31 | 0.13 | 0.32 | 0.16 | 0.80 | 0.23 | 0.73 | 0.25 | 0.79 | 0.21 |

Table 2
*Comparison of SCR, GCR and other methods for Example* 6.2

| | **SCR** | | **GCR** | | **SIR** | | **SAVE** | | **PHD** | |
|---|---|---|---|---|---|---|---|---|---|---|
| $\sigma$ | **DIST** | **SE** | **DIST** | **SE** | **DIST** | **SE** | **DIST** | **SE** | **DIST** | **SE** |
| 0.1 | 0.44 | 0.25 | 0.28 | 0.15 | 0.39 | 0.21 | 0.61 | 0.26 | 0.71 | 0.25 |
| 0.4 | 0.47 | 0.25 | 0.33 | 0.18 | 0.40 | 0.21 | 0.65 | 0.26 | 0.70 | 0.25 |
| 0.8 | 0.54 | 0.26 | 0.45 | 0.25 | 0.49 | 0.24 | 0.73 | 0.24 | 0.73 | 0.24 |

where $X$ and $\varepsilon$ are as defined in Example 6.1. Here, again, $q = 2$ and the central subspace is spanned by the vectors $(1, 0, 0, 0)^T$ and $(0, 1, 0, 0)^T$. We explore again the same grid of values for $\sigma$, using the same number of samples and sample size as in Example 6.1. Results are presented in Table 2.

We see that there is still a substantial improvement by GCR over SIR, SAVE and PHD. SIR slightly outperforms SCR, but the latter is much more accurate than SAVE and PHD.

In Section 5 we provided a population-level argument for the robustness of GCR against nonellipticity of the distribution of $X$. In the next example we compare GCR with OLS, PHD, SIR and SAVE when the distribution of $X$ is not elliptical.

EXAMPLE 6.3.  Consider the regression

$$Y = \sin^2(\pi X_2 + 1) + \sigma \varepsilon,$$

with predictor $X \in \mathbb{R}^4$ uniformly distributed on the set

$$[0, 1]^4 \setminus \{x \in \mathbb{R}^4 : x_i \leq 0.7, i = 1, 2, 3, 4\},$$

which defines a four-dimensional cube with a corner removed (this expedient is used to create an obvious asymmetry in the predictor distribution). We take again $\varepsilon \sim N(0, 1)$ and $\varepsilon \perp\!\!\!\perp X$. Here, the central subspace is of dimension $q = 1$ and is spanned by the vector $(0, 1, 0, 0)^T$. We perform the comparison once again drawing 500 samples of size $n = 100$ for each value $\sigma = 0.1$, 0.2 and 0.3. Results are presented in Table 3.

We see that GCR achieves a substantial improvement over OLS and PHD, and a modest one over SIR and SAVE. It also appears that SIR and SAVE are more robust than OLS and PHD against departures from ellipticity of $X$.

Next we compare contour regression methods with existing methods on instances where both the predictor dimension $p$ and the sample size $n$ are much larger than in the previous examples.

Example 6.4. We consider two cases with predictor dimension $p = 10$. We start with the regression

$$(18) \qquad Y = \cos(3X_1/2) + X_2^3/2 + \sigma\varepsilon,$$

where $X = (X_1, \ldots, X_{10})^T \sim N(0, I_{10})$, $\varepsilon \perp\!\!\!\perp X$ and $\varepsilon \sim N(0, 1)$. The central subspace has dimension $q = 2$ and is spanned by $(1, 0, \ldots, 0)^T$ and $(0, 1, \ldots, 0)^T$. As in Example 6.1, OLS is not considered in the comparison, as it can only detect one direction. The error standard deviation $\sigma$ is fixed at 0.1, 0.4 and 0.8 as in Examples 6.1 and 6.2, and for each such value we draw again 500 samples. Because of the increased dimension we now use a larger sample size; $n = 500$. The coefficients of the two terms $\cos(3X_1/2)$ and $X_2^3$ are chosen so that the "signal" for $X_2^3$ is strong for all $\sigma$ values, while the "signal" for $\cos(3X_1/2)$ is relatively weak for $\sigma = 0.8$. In this fashion, we can gather a sense of how the form of the regression function affects the performance of the various methods.

For both SCR and GCR we use $r = 5\%$ of the $500 \times 499/2 = 124{,}750$ empirical directions $(X_i - X_j)$ with the smallest $|Y_i - Y_j|$ or $\widehat{V}(\widehat{Z}_i, \widehat{Z}_j; \rho)$. For GCR the tube size is taken to be $\rho = 2$. Results are presented in Table 4.

For $\sigma = 0.1, 0.4$, both SCR and GCR, and especially GCR, achieve a marked improvement over the other methods. For $\sigma = 0.8$, though GCR still achieves some improvement, the accuracy of SCR is comparable to that of other estimators, indicating that under this level of noise the signal of

Table 3

*Comparison of GCR and other methods for Example 6.3*

| | GCR | | OLS | | PHD | | SIR | | SAVE | |
|---|---|---|---|---|---|---|---|---|---|---|
| $\sigma$ | DIST | SE | DIST | SE | DIST | SE | DIST | SE | DIST | SE |
| 0.1 | 0.10 | 0.05 | 0.17 | 0.07 | 0.24 | 0.10 | 0.13 | 0.06 | 0.14 | 0.08 |
| 0.2 | 0.12 | 0.06 | 0.19 | 0.09 | 0.29 | 0.12 | 0.18 | 0.08 | 0.22 | 0.12 |
| 0.3 | 0.20 | 0.14 | 0.22 | 0.10 | 0.36 | 0.16 | 0.22 | 0.10 | 0.34 | 0.20 |

Table 4

*Comparison of SCR, GCR and other methods for Example 6.4, regression (18)*

| | SCR | | GCR | | PHD | | SIR | | SAVE | |
|---|---|---|---|---|---|---|---|---|---|---|
| $\sigma$ | DIST | SE | DIST | SE | DIST | SE | DIST | SE | DIST | SE |
| 0.1 | 0.41 | 0.12 | 0.35 | 0.07 | 1.02 | 0.20 | 1.27 | 0.18 | 0.45 | 0.12 |
| 0.4 | 0.63 | 0.22 | 0.45 | 0.11 | 1.04 | 0.21 | 1.28 | 0.18 | 0.80 | 0.29 |
| 0.8 | 1.04 | 0.27 | 0.85 | 0.20 | 1.07 | 0.22 | 1.31 | 0.14 | 1.35 | 0.17 |

$\cos(3X_1/2)$ has dropped below the level detectable by most methods. We also observe that with sample size $n = 500$ the accuracy of SAVE has significantly increased compared with the previous examples where $n = 100$, suggesting that the relatively low accuracy of SAVE in the previous examples is probably due to its efficiency rather than the lack of population exhaustiveness.

The choice of $\rho = 2$ (compared to 1 used in the previous examples) is linked to the increased dimensionality: for large $p$ observations become sparse, and a thicker tube is needed to capture enough points. Experiments with numerous regression specifications invariably indicate that GCR with a relatively large tube size achieves outstanding improvements in accuracy. To benchmark the effect of $\rho$, the following simple quantity is useful: Let $X_1$, $X_2$ and $X_3$ be three independent observations from an $N(0, I_{10})$, and consider the probability of $X_3$ falling within the tube through $X_1$ and $X_2$ with $\rho = 2$. This probability is easily computed by simulation. From 500,000 simulated replicates of $X_1, X_2, X_3$ we find that

$$\Pr(d(X_3, \ell(X_1, X_2)) \leq 2) \approx 2.37\%.$$

Thus for sample size $n = 500$ there are on average $2 + 11.85 \approx 14$ observations in each tube. In contrast, if we take $\rho = 1$ the probability falls to approximately $9.4 \times 10^{-5}$, and the expected number of observations in each tube is 2.047, essentially equivalent to SCR.

Next, to confirm these results, we consider again the simpler regression in (17) and keep all the specifications of Example 6.1, except for taking $X \sim N(0, I_{10})$, $n = 500$, $r = 5\%$ and $\rho = 2$. Results are presented in Table 5.

We see that the broad patterns in Table 4 are confirmed, but the improvement by SCR and GCR appears to be more significant. The accuracy of SCR compared with GCR is increased somewhat, probably due to the simpler form of the regression.

In the next example we compare the performance of the above methods when the structural directions are determined by the variance, rather than the mean, function. This example demonstrates that, although the location structure (7) was postulated in deriving sufficient conditions for exhaustive estimation and robustness against nonellipticity, contour regression methods can be very effective also for regressions that are not based on location.

TABLE 5
*Comparison of SCR, GCR and other methods for Example 6.4, regression* (17)

|  | SCR | | GCR | | PHD | | SIR | | SAVE | |
|---|---|---|---|---|---|---|---|---|---|---|
| $\sigma$ | DIST | SE | DIST | SE | DIST | SE | DIST | SE | DIST | SE |
| 0.1 | 0.34 | 0.07 | 0.31 | 0.06 | 1.34 | 0.12 | 1.41 | 0.04 | 0.47 | 0.23 |
| 0.4 | 0.36 | 0.07 | 0.36 | 0.07 | 1.35 | 0.11 | 1.41 | 0.04 | 0.80 | 0.36 |
| 0.8 | 0.44 | 0.09 | 0.49 | 0.10 | 1.34 | 0.11 | 1.41 | 0.04 | 1.35 | 0.13 |

Table 6
*Comparison of SCR, GCR and other methods for Example 6.5*

| | SCR | | GCR | | PHD | | SIR | | SAVE | |
|---|---|---|---|---|---|---|---|---|---|---|
| $a$ | DIST | SE | DIST | SE | DIST | SE | DIST | SE | DIST | SE |
| 0.0 | 1.34 | 0.12 | 1.34 | 0.11 | 1.52 | 0.13 | 1.63 | 0.19 | 1.35 | 0.14 |
| 0.5 | 1.36 | 0.10 | 1.34 | 0.11 | 1.55 | 0.15 | 1.37 | 0.11 | 1.38 | 0.11 |
| 1.0 | 1.35 | 0.11 | 1.34 | 0.11 | 1.59 | 0.14 | 1.35 | 0.11 | 1.44 | 0.14 |

EXAMPLE 6.5. We consider the regression

$$Y = \tfrac{1}{2}(X_1 - a)^2 \varepsilon,$$

where $X \sim N(0, I_{10})$, $\varepsilon \sim N(0, 1)$ and $\varepsilon \perp\!\!\!\perp X$. Here the central subspace has dimension $q = 1$ and is spanned by $(1, 0, \ldots, 0)^T$. The variance of $Y$ is a quadratic function of $X_1$ centered at $a$, which is fixed at $a = 0, 0.5, 1$. Once again we generate 500 samples of size $n = 500$, and use $r = 5\%$ for SCR and GCR and $\rho = 2$ for GCR. Table 6 contains results for the comparison of SCR, GCR, PHD, SIR and SAVE. Cook and Li [5] proved that both OLS and PHD operate within the central *mean* subspace, and are therefore incapable of estimating a direction that only appears in the variance function. We include PHD in the comparison to serve as a benchmark for our subspace distance statistics.

Table 6 shows that contour regression methods are indeed capable of estimating the variance function direction, because their accuracy is much higher than the benchmark accuracy of PHD. Overall, the accuracy of SCR and GCR is similar to that of SIR and SAVE. We also observe that when $a$ is small, SAVE is more accurate than SIR, and the opposite is true when $a$ is large. The accuracy of contour methods does not appear to depend markedly on $a$. It is also worth mentioning that the errors in Table 6 are significantly larger than those in the previous examples regardless of the method used. This simply reflects the fact that estimating variance structures is more difficult than estimating mean structures.

Finally, we come to the issue of estimating the structural dimension $q$. As mentioned in Section 2.4, we believe that an asymptotic test for SCR can be developed along the same lines employed by existing methods such as SIR and PHD, though this will require work beyond the scope of the present paper. The development of an asymptotic test for GCR would hinge on an asymptotic analysis of the GCR estimator, which has not been pursued here. Nevertheless, we can empirically assess the capability of SCR and GCR to estimate $q$ by examining how much the eigenvalues corresponding to the central subspace are separated from those corresponding to its complement (i.e., the contour space).

For SCR we use the matrix

$$2I_p - \widehat{\Sigma}^{-1/2} \widehat{K}(c) \widehat{\Sigma}^{-1/2}. \tag{19}$$

This is the sample version of the population matrix $2I_p - \Sigma^{-1/2} K(c) \Sigma^{-1/2}$. From Theorem 2.3 we know that the eigenvalues of the population matrix corresponding to contour directions are identically 0 and those corresponding to the central subspace are strictly positive. Thus we expect the eigenvalues of the sample matrix (19) to behave similarly. We consider again the simulations for regressions (18) and (17) in Example 6.4, with $\sigma = 0.4$. We compute the ten eigenvalues of the matrix (19) for each of the 500 samples, say $\hat{\lambda}_{\ell,1}, \ldots, \hat{\lambda}_{\ell,10}$. From each of the ten sets of simulated eigenvalues $\{\hat{\lambda}_{1,j}, \ldots, \hat{\lambda}_{500,j}\}$, $j = 1, \ldots, 10$, we then compute an average $\hat{\lambda}_{\cdot j}$ and a standard error $\hat{\tau}_j$. These numbers are shown in the two SCR columns of Table 7.

For GCR we use the matrix

$$\widehat{G}(c) = \sum_{(i,j) \in N} (\widehat{Z}_i - \widehat{Z}_j)(\widehat{Z}_i - \widehat{Z}_j)^T$$
$$\times I(\widehat{V}(\widehat{Z}_i, \widehat{Z}_j; \rho) \leq c) \Big/ \sum_{(i,j) \in N} I(\widehat{V}(\widehat{Z}_i, \widehat{Z}_j; \rho) \leq c),$$

where the index set $N$ is as defined in (4). This matrix is proportional to $\widehat{F}(c)$ in (12), rescaled so that it estimates $E[(Z_i - Z_j)(Z_i - Z_j)^T \mid |I_i - I_j| \leq c]$ instead of $E[(Z_i - Z_j)(Z_i - Z_j)^T I(|I_i - I_j| \leq c)]$. Though we did not prove a theorem for GCR analogous to Theorem 2.3, we mimic the SCR case and compare the eigenvalues of $2I_p - \widehat{G}(c)$ (it is the separation of eigenvalues that matters here). The simulation averages and standard errors of these eigenvalues over 500 samples for regressions (18) and (17) (with $\sigma = 0.4$) are shown in the two GCR columns of Table 7.

From Table 7 we see that for both methods, and in both regressions, the eigenvalues $\hat{\lambda}_{9,j}$ and $\hat{\lambda}_{10,j}$, which correspond to vectors in the central subspace, are significantly larger than the other eigenvalues. Furthermore, the contrast between $\hat{\lambda}_{9,j}$ and $\hat{\lambda}_{10,j}$ and the remaining eigenvalues appears to be stronger for GCR than for SCR, suggesting that GCR is more sensitive in identifying the central subspace.

**7. An application.** We consider a data set concerning the effect on soil evaporation of various air and soil conditions such as temperature, humidity and wind speed ([12]; it is available in the *Arc* package—see http://www.stat.umn.edu/arc/software.html). There are $p = 10$ predictors: average daily air temperature (Avat), area under the daily humidity curve (Avh), area under the daily soil temperature curve (Avst), maximum daily air temperature (Maxat), maximum daily humidity (Maxh), maximum daily soil temperature (Maxst), minimum daily air temperature (Minat), minimum daily humidity (Minh), minimum daily soil temperature (Minst) and total wind speed in miles/hour (Wind). The response is daily soil evaporation (Evap). The data are collected over a period of 46 days, but do not show any obvious serial dependence. Hence for simplicity we treat the data as independent replicates with $n = 46$.

Figure 2 is the scatterplot matrix of the ten predictors, which does not seem to suggest serious departures from ellipticity. Furthermore, simultaneous Box–Cox transformations of these predictors do not lead to significant improvements in ellipticity. Hence we use the untransformed predictors for our analysis. We apply both SIR and GCR to the data, using the negative Evap as $Y$. The two upper panels of Figure 3 are the scatterplots of $Y$ versus the first two SIR directions, SIR1 and SIR2, on the standardized scale $\widehat{Z}$. The scatterplot for $Y$ versus SIR1 (upper-left panel) shows a strong monotone trend which is almost linear. In contrast, the scatterplot of $Y$ versus SIR2 (upper-right panel) does not show a detectable pattern. The two lower

TABLE 7
*Averages (EVAL) and standard errors (SE) of eigenvalues from SCR and GCR*

| EVAL(SE) | MODEL I | | MODEL II | |
|---|---|---|---|---|
| | **SCR** | **GCR** | **SCR** | **GCR** |
| $\hat{\lambda}_{\cdot 1}\ (\hat{\tau}_1)$ | $-0.26\ (0.05)$ | $-0.55\ (0.12)$ | $-0.23\ (0.04)$ | $-0.48\ (0.11)$ |
| $\hat{\lambda}_{\cdot 2}\ (\hat{\tau}_2)$ | $-0.18\ (0.04)$ | $-0.37\ (0.10)$ | $-0.15\ (0.04)$ | $-0.32\ (0.09)$ |
| $\hat{\lambda}_{\cdot 3}\ (\hat{\tau}_3)$ | $-0.11\ (0.04)$ | $-0.23\ (0.08)$ | $-0.09\ (0.03)$ | $-0.21\ (0.08)$ |
| $\hat{\lambda}_{\cdot 4}\ (\hat{\tau}_4)$ | $-0.05\ (0.04)$ | $-0.11\ (0.08)$ | $-0.04\ (0.03)$ | $-0.10\ (0.07)$ |
| $\hat{\lambda}_{\cdot 5}\ (\hat{\tau}_5)$ | $0.01\ (0.04)$ | $0.00\ (0.07)$ | $0.02\ (0.04)$ | $0.00\ (0.06)$ |
| $\hat{\lambda}_{\cdot 6}\ (\hat{\tau}_6)$ | $0.07\ (0.04)$ | $0.11\ (0.07)$ | $0.07\ (0.04)$ | $0.10\ (0.07)$ |
| $\hat{\lambda}_{\cdot 7}\ (\hat{\tau}_7)$ | $0.14\ (0.04)$ | $0.23\ (0.08)$ | $0.13\ (0.04)$ | $0.20\ (0.07)$ |
| $\hat{\lambda}_{\cdot 8}\ (\hat{\tau}_8)$ | $0.23\ (0.05)$ | $0.37\ (0.08)$ | $0.21\ (0.05)$ | $0.33\ (0.07)$ |
| $\hat{\lambda}_{\cdot 9}\ (\hat{\tau}_9)$ | $0.41\ (0.08)$ | $0.91\ (0.11)$ | $0.72\ (0.07)$ | $1.08\ (0.09)$ |
| $\hat{\lambda}_{\cdot 10}\ (\hat{\tau}_{10})$ | $1.17\ (0.06)$ | $1.23\ (0.07)$ | $1.14\ (0.07)$ | $1.21\ (0.07)$ |

MODEL I is regression (18) with $\sigma = 0.4$, and MODEL II is regression (17) with $\sigma = 0.4$ and a ten-dimensional predictor $X$.

panels of Figure 3 are the scatterplots of $Y$ versus the first two GCR directions, GCR1 and GCR2, on the standardized scale $\widehat{Z}$. The plot for $Y$ versus GCR1 (lower-left panel) also shows a clear monotone, but slightly nonlinear, trend. What is interesting, however, is that the scatterplot of $Y$ versus GCR2 (lower-right panel) suggests a $U$-shaped pattern. A three-dimensional spin plot of $Y$ versus (GCR1, GCR2) shows a mean surface that, roughly speaking, is folded in the GCR2 direction and tilted upwards in the GCR1 direction. In the $Y$-versus-GCR2 scatterplot, five points (labeled by "+") sit above the $U$-shape near GCR2 $= 0$, appearing to weaken the $U$-shaped pattern. However, these points are far out in the direction of GCR1 with high values of $Y$—corresponding to the five points labeled by "+" in the perspective scatterplot in Figure 4.

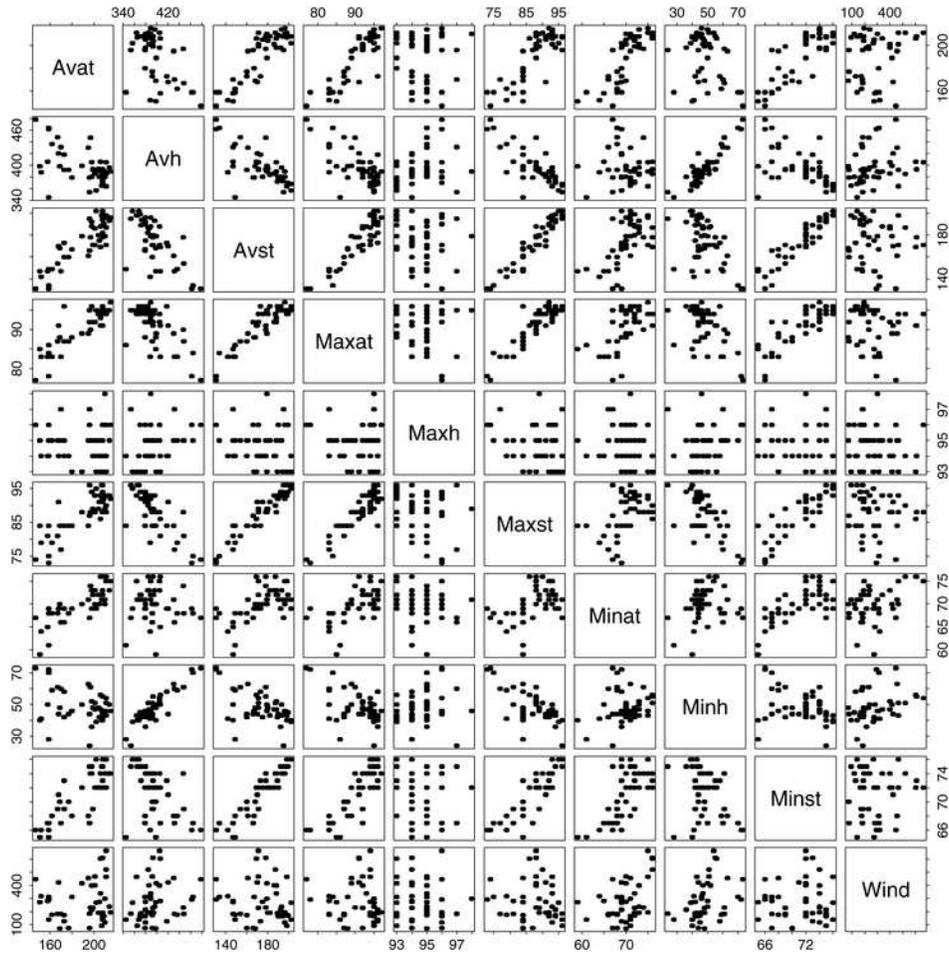

Fig. 2. *Scatterplot matrix for the predictors in the soil evaporation data.*

Because in this data set $p = 10$ and $n = 46$, as discussed in Section 6 (following Example 6.4) we need to choose a rather large radius $\rho$ to capture enough points in each tube. For this application of GCR we used $\rho = 3.5$ (on the $\widehat{Z}$-scale) and 15% of the $\binom{46}{2} = 1035$ pairs (i.e., 155 pairs) of points among which $\widehat{V}(\widehat{Z}_i, \widehat{Z}_j; \rho = 3.5)$ are the smallest. For SIR we used six slices defined so as to contain roughly the same number of points.

Although without a formal testing procedure we cannot yet determine the statistical significance of GCR2, the 2D and 3D scatterplots from our GCR analysis do suggest that a second direction might be relevant in the evaporation data. Due to its small sample size relative to the dimension of the predictor, this example does not allow us to draw strong conclusions, but it demonstrates once again that GCR is more sensitive than classical meth-ods in detecting complex regression surfaces—in this instance monotone in one direction and $U$-shaped in another. This was anticipated by theoretical analysis in Section 4 and supported by simulation studies in Section 6.

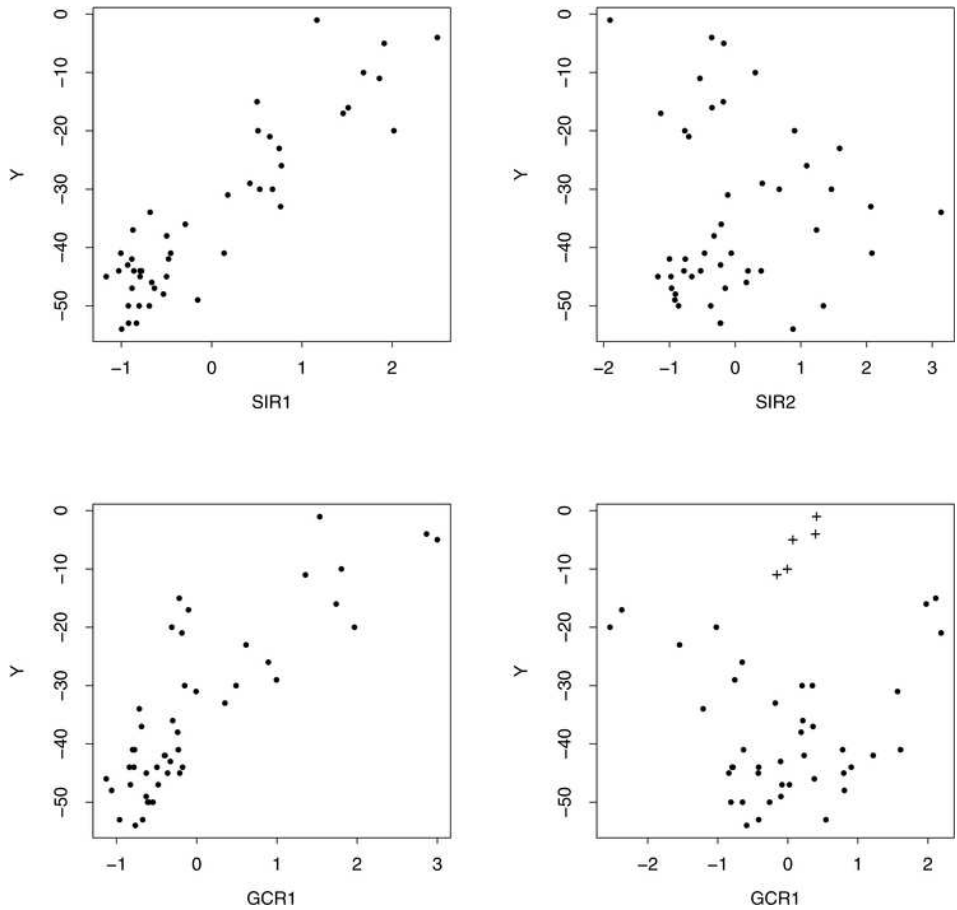

FIG. 3. *Scatterplots of the response* (−Evap) *versus the first two SIR directions* (upper panels) *and the first two GCR directions* (lower panels) *for the soil evaporation data.*

**8. Conclusions.** The contour regression methods introduced in this paper have strength in several aspects. First, under mild conditions they achieve exhaustive estimation of the central subspace at the $\sqrt{n}$-convergence rate. In comparison with existing global estimators such as OLS, PHD and SIR, contour regression estimators are more comprehensive, capable of picking up all directions in the central subspace without relying on special response patterns (e.g., monotone or $U$-shaped trends). In particular, GCR achieves exhaustiveness essentially without any assumption other than that of ellipticity of the distribution of $X$. Second, by design contour regression methods are capable of exploiting interslice information which is not accessible to methods based on slicing. This partly explains their improved accuracy over SIR and SAVE, which we have discussed from an analytical standpoint and documented through simulations. In fact, we think that the advantage of contour methods over SAVE is due to the gain in efficiency such as achieved via the use of interslice information, and not the structural inability of SAVE to capture linear trends. Third, GCR achieves a degree of robustness against nonellipticity of the distribution of $X$. In this respect contour regression is akin to the adaptive methods mentioned in the Introduction. Unlike adaptive methods, however, contour methods are computationally simple, the level of computational burden being essentially that of principal component analysis. In particular, they do not require iterative maximization of a multivariate nonparametric function, which can be a substantial advantage, especially if the dimension is large, or if multiple local maxima are present in the iterative maximization.

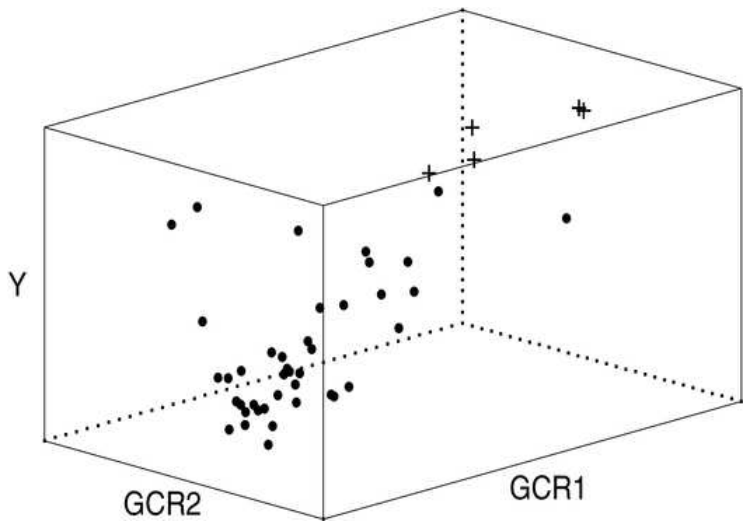

FIG. 4.    *A view of the* 3*D plot of the response (*−Evap*) versus the first two GCR directions in the soil evaporation data.*

Because contour vectors are extracted according to a threshold on response variation, our methods are logically analogous to a one-dimensional kernel or nearest-neighbor estimator. If the distribution of the predictor $X$ is elliptical, the threshold need not go to zero in our asymptotic arguments for SCR, which makes it possible to achieve the $\sqrt{n}$-rate regardless of the dimensions $p$ and $q$. In this respect contour regression is similar to traditional global methods such as OLS, PHD, SIR and SAVE. However, if ellipticity fails and/or is difficult to establish through predictor transformations, we can employ a relatively small threshold, operating in a spirit more similar to that of adaptive methods.

We do not claim that contour regression estimators will outperform other methods under all circumstances. For example, OLS is the maximum likelihood estimator if the regression surface is linear and the error term is normal, and tends to perform very well if the surface is nearly linear or clearly monotone. Similarly favorable circumstances exist for PHD, SIR and SAVE as well.

The ideas of contour regression raise many questions that have not been addressed within the the scope of this paper. In particular, the asymptotic properties of GCR, as well as test statistics for estimating the structural dimension $q$, have not yet been developed. We do expect that $\sqrt{n}$-convergence can be achieved by GCR if the threshold $c$ is taken as fixed, because this in effect includes in the computation a number of line segments proportional to the total number of observation pairs. We also expect that test statistics for determining $q$ can be constructed based on Theorem 2.3, along lines similar to those in [2, 14]. Also, we have not provided a systematic method for choosing the thresholding constant $c$ (or the ratio $r$) for SCR and GCR, as well as the tube radius $\rho$ for GCR, which should ideally be based on data-driven criteria. Other useful developments will concern the asymptotic behavior of GCR when the threshold $c$ is allowed to go to zero as the sample size $n$ tends to infinity. Theorem 5.1 suggests that even without ellipticity of $X$ the correct asymptotic behavior would still be guaranteed. However, in this case we do not expect a $\sqrt{n}$-convergence rate—at least not for all structural dimensions. To further improve efficiency it may be helpful to experiment with windows other than the current rectangular ones in selecting contour vectors. It may also be possible to apply an idea similar to local linear regression [11] to correct the possible edge effect caused by the line segments lying in the outskirts of the data cloud. Another worthwhile line of research would be to extend contour methods to dependent data, for example, to weakly dependent Gaussian time series (see [23]). Finally, as we have seen from Example 2.3 (binary response), contour methods do apply to discrete numerical responses. Moving forward along this direction, we could generalize contour methods to ordinal categorical or even purely categorical responses. This will require appropriate renditions for the concepts of *absolute difference* (e.g., "the same" and "not the same") and *distributional spread* (e.g., a concentration index) for categorical modalities; these could be used to define contour direction identifiers for SCR and GCR, respectively. We leave these issues to future studies.

## APPENDIX

PROOF OF LEMMA 2.1. (a) Let $f_{ij}$ denote the joint densities of $(V_i, V_j)$, and so on. For example, $f_1$ is the density of $V_1$ and $f_{123}$ is the joint density of $(V_1, V_2, V_3)$. Similarly, let $f_{ij|k}$ and so on denote conditional densities. For example, $f_{23|4}$ is the conditional joint density of $(V_2, V_3)$ given $V_4$. We need to show that

$$(20) \quad f_{1245|36}(v_1, v_2, v_4, v_5|v_3, v_6) = f_{14|36}(v_1, v_4|v_3, v_6)f_{25|36}(v_2, v_5|v_3, v_6).$$

Without loss of clarity we can omit $v_1, \ldots, v_6$ from the density. Thus the above equality becomes $f_{1245|36} = f_{14|36}f_{25|36}$. The left-hand side of (20) is

$$f_{123456}/f_{36}.$$

Because $(V_1, V_2, V_3) \perp\!\!\!\perp (V_4, V_5, V_6)$, $V_1 \perp\!\!\!\perp V_2|V_3$ and $V_4 \perp\!\!\!\perp V_5|V_6$, the numerator in the above ratio is factorized into $f_{1|3}f_{2|3}f_3f_{4|6}f_{5|6}f_6 = f_{13}f_{23}f_{46}f_{56}/(f_3f_6)$, and the denominator is factorized into $f_3f_6$. Thus the left-hand side reduces to $f_{13}f_{23}f_{46}f_{56}/(f_3f_6)^2$. The right-hand side of (20) is the ratio $f_{1436}f_{2536}/f_{36}^2$. Because $(V_1, V_2, V_3) \perp\!\!\!\perp (V_4, V_5, V_6)$ this ratio becomes $f_{13}f_{46}f_{23}f_{56}/(f_3f_6)^2$, the same quantity to which the left-hand side of (20) is reduced.

(b) Suppose $(V_1, V_2) \perp\!\!\!\perp (V_3, V_4)$. We want to show that

$$(21) \quad\quad\quad\quad\quad\quad f_{13|24} = f_{1|24}f_{3|24}.$$

The left-hand side is $f_{1234}/f_{24}$, which, because $(V_1, V_2) \perp\!\!\!\perp (V_3, V_4)$, reduces to $f_{12}f_{34}/(f_2f_4)$. The right-hand side of (21) is $f_{124}f_{324}/(f_{24})^2 = f_{12}f_4f_2f_{34}/(f_2f_4)^2 = f_{12}f_{34}/(f_2f_4)$, which completes the proof. □

PROOF OF THEOREM 4.1. Because the proof is basically the same as that of Theorem 2.1, we only highlight the differences. There is no change in paragraphs 1, 2, 3, 4, 6 of the proof of Theorem 2.1 except for replacing, wherever applicable, $|\tilde{Y} - Y| \le c$ by $V(X, \tilde{X}) \le c$, $K(c)$ by $G(c)$, $K_1(c)$ by

$$G_1(c) = E[(Z - \tilde{Z})(Z - \tilde{Z})^T|V(\tilde{Z}, Z) \le c],$$

and "Assumption 2.1" by "Assumption 4.1."

Replace the fifth paragraph by the following argument: Because $(Z, Y)$ and $(\phi_i(Z), Y)$ have the same distribution, and because $(Z, Y)$ and $(\tilde{Z}, \tilde{Y})$ are independent, the distributions of $(Z, Y, \tilde{Z}, \tilde{Y})$ and $(\phi_i(Z), Y, \phi_i(\tilde{Z}), \tilde{Y})$ are identical. Hence

$$E[(\tilde{Z} - Z)(\tilde{Z} - Z)^T | V(Z, \tilde{Z}) \le c]$$
$$= E[(\phi_i(\tilde{Z}) - \phi_i(Z))(\phi_i(\tilde{Z}) - \phi_i(Z))^T | V(\phi_i(Z), \phi_i(\tilde{Z})) \le c].$$

We claim that, for any $a$ and $b$ in $\mathbb{R}^p$, $V(a, b) = V(\phi_i(a), \phi_i(b))$. By definition,

$$V(a, b) = \text{var}(Y | Z = (1 - t)a + tb \text{ for some } t).$$

Because $(Z, Y)$ and $(\phi_i(Z), Y)$ have the same distribution, the condition in the above conditional variance can be replaced by $\phi_i(Z) = (1 - t)a + tb$ or $Z = \phi_i^{-1}((1 - t)a + tb)$. Because, as we have noted, $\phi_i = \phi_i^{-1}$, and also because $\phi_i: \mathbb{R}^p \mapsto \mathbb{R}^p$ is a linear function, the condition can be replaced by $Z = (1 - t)\phi_i(a) + t\phi_i(b)$. Therefore,

$$V(a, b) = \text{var}(Y | Z = (1 - t)\phi_i(a) + t\phi_i(b) \text{ for some } t) = V(\phi_i(a), \phi_i(b)),$$

as claimed. Consequently,

$$E[(\tilde{Z} - Z)(\tilde{Z} - Z)^T | V(Z, \tilde{Z}) \le c]$$
$$= E[(\phi_i(\tilde{Z}) - \phi_i(Z))(\phi_i(\tilde{Z}) - \phi_i(Z))^T | V(Z, \tilde{Z}) \le c].$$

Now follow through the rest of the fifth paragraph in the proof of Theorem 2.1, replacing $|Y - \tilde{Y}| \le c$ by $V(X, \tilde{X}) \le c$ in one place. $\square$

PROOF OF LEMMA 3.2. By Fubini's theorem we have

$$E(S) - E(T) = \int_0^\infty \int_t^\infty \left( \frac{p(s)}{q(s)} - 1 \right) q(s) \, ds \, dt \equiv \int_0^\infty G(t) \, dt.$$

We now show that $G(t) < 0$ for all $t > 0$. Because $p(s)/q(s)$ is a decreasing function and because

$$\int_0^\infty \left( \frac{p(s)}{q(s)} \right) q(s) \, ds = 1,$$

the function $p(s)/q(s)$ is greater than 1 at $s = 0$, equal to 1 at some $s_0 > 0$ and less than 1 afterward. Hence $G'(t) = q(t) - p(t)$ is less than 0 for $t < s_0$ and greater than 0 for $t > s_0$. So $G(t)$ first decreases and then increases. However, it is easy to see that $G(0) = 0$ and $\lim_{t \to \infty} G(t) = 0$. Hence $G(t) < 0$ for all $t > 0$. $\square$

LEMMA A.1. *Let $\mathcal{S}_1$ and $\mathcal{S}_2$ be two linear subspaces of $\mathbb{R}^p$ and let $A$ be an open set in $\mathbb{R}^p$ containing the origin. Suppose $A \cap \mathcal{S}_1 \subset A \cap \mathcal{S}_2$. Then $\mathcal{S}_1 \subset \mathcal{S}_2$.*

PROOF.    Let $v$ be a vector in $\mathcal{S}_1$. Because $A$ is an open set containing the origin, for sufficiently small $\lambda > 0$, $\lambda v \in A$. Hence $\lambda v \in A \cap \mathcal{S}_1$. Because $A \cap \mathcal{S}_1 \subset A \cap \mathcal{S}_2$, $\lambda v$ also belongs to $A \cap \mathcal{S}_2$. Hence $\lambda v$ belongs to $\mathcal{S}_2$. Because $\mathcal{S}_2$ is a linear subspace, $v$ belongs to $\mathcal{S}_2$.    $\square$

**Acknowledgments.**    We are grateful to two referees and an Associate Editor for their insightful and thorough reviews, which led to significant improvements to this article.

... 

B. Li
F. Chiaromonte
Department of Statistics
326 Thomas Building
Pennsylvania State University
University Park, Pennsylvania 16802-2111
USA
e-mail: bing@stat.psu.edu
e-mail: chiaro@stat.psu.edu

H. Zha
Department of Computer
  Science and Engineering
343F IST Building
Pennsylvania State University
University Park, Pennsylvania 16802
USA
e-mail: zha@cse.psu.edu